\documentclass[a4paper,article,10pt]{memoir}
\usepackage[english]{babel}
\usepackage{CCLAuthor}
\usepackage[reqno]{amsmath}
\usepackage{amssymb}
\usepackage{amsthm}
\usepackage{amsmath}
\usepackage{mathtools}
\usepackage{color}
\usepackage{multirow}
\usepackage{url}
\usepackage{algorithm2e}
\usepackage{hyperref} 
\usepackage[figuresright]{rotating}
\theoremstyle{plain}
\newtheorem{theorem}{Theorem}

\theoremstyle{definition}

\newtheorem{example}[theorem]{Example}
\numberwithin{theorem}{chapter}
\DeclareMathOperator*{\argmin}{\arg\!\min}
\usepackage[round,authoryear]{natbib} 
 \usepackage[]{graphicx}              
\usepackage{rotating}
\usepackage{adjustbox}
\usepackage{float}
\usepackage{titling}
\thanksmarkseries{alph}
\begin{document}
%
%
%
\title{Data-driven Methods for Delay Differential Equations}
%
%
\author{%
    Dimitri Breda\textsuperscript{*}\thanks{Dimitri Breda is member of INdAM Research group GNCS and of UMI Research group ``Mo\-del\-li\-sti\-ca socio-epidemiologica''. His work was partially supported by the Italian Ministry of University and Research (MUR) through the PRIN 2022 project (No. 20229P2HEA) ``Stochastic numerical modelling for sustainable innovation'' (CUP: G53C24000710006).}, Xunbi A. Ji\textsuperscript{\dag}\thanks{Xunbi Ji was supported by a  Beyster Computational Innovation Graduate Fellowship and by a Rackham Graduate Student Research Grant.},
    G\'abor Orosz\textsuperscript{\dag\ddag}\thanks{G\'abor Orosz was supported by a Fulbright Fellowship.}, 
    and Muhammad Tanveer\textsuperscript{*}\thanks{Muhammad Tanveer is member of INdAM Research group GNCS and his work was supported by the Italian Ministry of University and Research (MUR) through a PhD grant PNRR DM351/22.}
    \\ \smallskip\small
   \textsuperscript{*} 
    CDLab – Computational Dynamics Laboratory,\\ Department of Mathematics, Computer Science and Physics,\\ University of Udine, Udine, Italy
    \\ 
    \textsuperscript{\dag} 
    Department of Mechanical Engineering, University of Michigan, Ann Arbor, USA
    \\ 
    \textsuperscript{\ddag} 
    Department of Civil and Environmental Engineering, University of Michigan, Ann Arbor, USA
    }
   
%
    \maketitle
%
%
%
\begin{abstract}
Data-driven methodologies are nowadays ubiquitous. 
Their rapid development and spread have led to applications even beyond the traditional fields of science. 
As far as dynamical systems and differential equations are concerned, neural networks and sparse identification tools have emerged as powerful approaches to recover the governing equations from available temporal data series. 
In this chapter we first  illustrate possible extensions of the sparse identification of nonlinear dynamics (SINDy) algorithm, originally developed for ordinary differential equations (ODEs), to delay differential equations (DDEs) with discrete, possibly multiple and unknown delays. 
Two methods are presented for SINDy, one directly tackles the underlying DDE and the other acts on the system of ODEs approximating the DDE through pseudospectral collocation. 
We also introduce another way of capturing the dynamics of DDEs using neural networks and trainable delays in continuous time, and present the training algorithms developed for these neural delay differential equations (NDDEs).
The relevant MATLAB implementations for both the SINDy approach and for the NDDE approach are provided.
These approaches are tested on several examples, including classical systems such as the delay logistic and the Mackey-Glass equation, and directly compared to each other on the delayed R\"ossler system.
We provide insights on the connection between the approaches and future directions on developing data-driven methods for time delay systems.
\end{abstract}


\CCLsection{Introduction}

Data-driven techniques have revolutionized model development for complex systems, given the availability of large-scale datasets and powerful computational resources \citep{bpk16}. 
This revolution is especially applicable in domains where interactions are not only nonlinear and high-dimensional but also inherently delayed \citep{hale77,diekmann95}, thus exhibiting memory effects. 
In these cases, traditional modeling approaches based on first-principles formulations often fail to capture the complicated temporal dependencies that are clearly manifested in the data. 
This generates an increasing need for new data-driven methodologies which directly derive (rather than assume) the governing equations with enhanced predictive capability from observational data, and provide insight into the complex mechanisms involved. 

While very effective in identifying governing equations for ordinary differential equations (ODEs), criticism has continued to exist in the case of systems exhibiting time-delayed dynamics \citep{smith2011introduction}. 
In many real applications, the evolution of the system depends both on its present and past states, which can be captured by delay differential equations (DDEs) \citep{hale77,diekmann95}. 
However, the intrinsic nature of the delay---often manifesting itself as unobservable latencies due to signal propagation, feedback loops, or processes themselves \citep{erneux2009applied}---renders the task of reconstruction more challenging. Traditional approaches tend to focus on knowledge predetermined about these delays in diverse scenarios such as physiological processes \citep{Hutchinson1948}, neural activity \citep{shayer2000} and climate dynamics \citep{saunders2001boolean}. 
This highlights the need for robust data-driven methods, capable of simultaneously identifying the dynamic behavior and the unknown delay parameters from empirical data without relying on prior knowledge. 
Methods of this kind are needed not only to enhance predictive performance, but also to deepen our understanding toward complex long-term interactions characterizing many natural and engineered systems \citep{stepan1989retarded,insperger2011semi}.

Among many data-driven methods, two complementary approaches have recently emerged for dealing with delayed dynamics in real-world systems, each exploiting the advantages of data-driven methods in discovering governing equations. 
The first is based on sparse regression and it is called sparse identification of nonlinear dynamics (SINDy), as discussed by \cite{bpk16}. 
This method generates interpretable models through the selection of a handful of active terms from a collection of candidate functions. 
Initially developed for ordinary differential equations (ODEs) \citep{bpk16}, SINDy has since been extended to more general classes of problems, including partial differential equations \citep{rud17}, stochastic differential equations \citep{bnc18} and DDEs with discrete delays \citep{bbt24,kopeczi23,pec24,sandoz23}. Recent extensions \citep{bbt24}, named Expert SINDy ({E-SINDy}) and Pragmatic SINDy ({P-SINDy}), cope with DDEs either by incorporating the known delay structures into the library or by transforming the delayed system into an ordinary differential system via pseudospectral collocation \citep{bdgsv16, bdls16}. 
Parallel to these approaches, neural network-based methods arose~\citep{pei2013mapping,bk19}, representing flexible gradient-based solutions that can simultaneously extract both the dynamics and the delay parameters from raw data. 
In particular, neural delay differential equations (NDDEs) use architectures with trainable delays that are particularly useful in capturing complicated temporal dependencies arising in systems with memory effects~\citep{ji2020,gupta2023neural,levine2022framework}.

This chapter aims to provide a hands-on, comprehensive overview of the above-mentioned data-driven approaches for DDEs from time-series data. 
In doing so, it endeavors to take an entirely different view compared to the deterministic, first-principle-based methods that dominate the rest of this volume. 
We exhibit and compare SINDy and NDDEs for reconstructing and simulating delayed dynamics, while achieving some compromise among interpretability, computational efficiency and predictive accuracy.
It is hoped that such an integration will lead to robust, interpretable and computationally efficient tools for capturing rich memory-dependent dynamics often exhibited by real-world systems in motion. 
Apart from laying down the theoretical foundations, MATLAB code examples and implementation guidelines, designed to be accessible by a general audience, are given throughout the chapter. 
The chapter aspires to bridge the gap between classical deterministic modeling of time delay systems and the burgeoning realm of data-driven discovery, thereby fostering a deeper understanding of complex dynamical phenomena in different fields.

This chapter is organized to provide the reader with a walk-through of the topic. 
Section~\ref{sec:2} defines the fundamental properties of DDEs establishing key concepts followed by the sparse identification of DDEs.
The latter includes {E-SINDy} and {P-SINDy}, each accompanied by the examples. 
In Section~\ref{sec:NN}, we dive into the neural network techniques for delayed systems and present the formulation of NDDEs with trainable delays including details on loss function design and the training methods. 
Practical implementation is emphasized throughout the chapter and links to MATLAB code examples are provided to ensure that the material is accessible to a broad audience. 
Finally, Section~\ref{sec:compare} presents comparative analyses and discussions to highlight the trade-offs between interpretability and computational efficiency, concluding with insights into future research directions in the data-driven discovery of dynamical systems with memory.

\CCLsection{Sparse Identification of Delay Equations}\label{sec:2}

We begin in Section~\ref{sec:DDEs} with a brief overview of DDEs, discussing the basic mathematical formulation and key properties that we need in the following treatment. 
Next, in Section~\ref{sec:sindy}, we introduce the SINDy framework for ODEs, that it will be then extended to DDEs. 
In Section~\ref{sec:Etauknown}, we explore E-SINDy, the ``expert'' approach to SINDy for DDEs when delays are known, followed by an external optimization methods for identifying unknown delays and non-multiplicative parameters in Section~\ref{sec:Etauunknown}. 
We then introduce in Section~\ref{sec:PSINDy} P-SINDy, the ``pragmatic'' approach, a more computationally efficient alternative. 
Basic numerical experiments are provided at the end of every section to illustrate the introduced methodologies. 
In Section~\ref{sec:discussion}, we provide more numerical experiments to explore the performance of the different techniques for identifying and reconstructing DDEs from data.

\CCLsubsection{Delay Differential Equations}\label{sec:DDEs}

Here we summarize the basic aspects of initial value problems (IVPs) for delay systems in order to introduce their description in terms of dynamical systems. 
In this respect, we recall the main ingredients of the relevant semigroup theory from \citep{engel2000one}, see also Chapter 3 of \citep{breda2022controlling}.
For a full theoretical treatment of DDEs see, e.g., \citep{hale77, diekmann95}. 

For ${\tau > 0}$ real and $n$ a positive integer, let ${\mathcal{X} \coloneqq C([-\tau, 0], \mathbb{R}^n)}$
be the space of continuous functions ${[-\tau, 0]\to\mathbb{R}^n}$, which is a Banach space when equipped with the uniform norm
${\| \varphi \|_\mathcal{X} \coloneqq \max_{s \in [-\tau, 0]} \| \varphi(s) \|}$ where ${\| \cdot \|}$ is any norm in $\mathbb{R}^n$. 
A general nonlinear autonomous DDE takes the form
\begin{equation}\label{eq:ddef}
x'(t) = f(x_t),
\end{equation}
where ${x\in \mathbb{R}^n}$, prime denotes the (right hand) derivative with respect to time $t$, the right hand side (RHS) ${f: \mathcal{X} \to \mathbb{R}^n}$ governs the system's dynamics and the \textit{history function} ${x_t \in \mathcal{X}}$ at time $t$ reads ${x_t(s) \coloneqq x(t + s)}$, ${s \in [-\tau, 0]}$, following the standard notation \citep{krasovskii63}.
Here, $\tau$ represents the maximum delay, meaning that the system's current evolution depends on its past over the interval ${[-\tau, 0]}$. 
An IVP for \eqref{eq:ddef} is specified as
\begin{equation}\label{eq:IVPDDEF}
\left\{\setlength\arraycolsep{0.1em}\begin{array}{ll}
x'(t)=f(x_t),&\quad t \geq 0,\\[1mm]
x(s)=\varphi(s),&\quad s\in[-\tau,0],
\end{array}\right.
\end{equation}
for ${\varphi\in\mathcal{X}}$ playing the role of the initial history, viz.\ the state $x_{0}$ at time 0. 
The RHS $f$ is assumed to be sufficiently smooth to ensure that \eqref{eq:IVPDDEF} has a unique solution depending smoothly on $\varphi$. 
For general well-posedness results see, e.g., \citep{bellen2009recent}.

It is clear from \eqref{eq:IVPDDEF} that DDEs generate infinite-dimensional dynamical systems on the state space $\mathcal{X}$. 
Due to this infinite-dimensional nature, DDEs can exhibit complex and rich dynamics already in the scalar (${n=1}$) case \citep{spjc07,mg77}. 
In the following, we refer to $n$ as to the \textit{physical dimension} (i.e., the number of equations) and to the infinite dimension of the state space $\mathcal{X}$ as to the \textit{dynamical dimension} (i.e., the classical ``amount'' of information needed to progress with the time evolution). 
Next we describe the dynamical systems framework by recalling the main tools from semigroup theory.

A well-posed IVP \eqref{eq:IVPDDEF} allows one to introduce the solution operator ${T(t) :  \mathcal{X} \to  \mathcal{X}}$ at time ${t\geq0}$ as ${T(t)\varphi = x_t}$,
which associates to the initial state $\varphi$ the corresponding history at time $t$. 
The operator $T(t)$ acts as a translation operator, with extension rule defined through \eqref{eq:ddef}. 
The collection ${\{T(t)\}_{t \geq 0}}$ forms a strongly continuous semigroup of bounded operators. 
Consequently, we can define the relevant infinitesimal generator ${\mathcal{A} : \mathcal{D}(\mathcal{A})\subseteq\mathcal{X}\to\mathcal{X}}$ as the unbounded operator given by
\begin{equation}\label{eq:generator}
    \mathcal{A}\varphi=\varphi',\qquad\mathcal{D}(\mathcal{A})=\{\varphi\in\mathcal{X}:\varphi'\in\mathcal{X}\text{ and }\varphi'(0)=f(\varphi)\}.
\end{equation}
Note that the underlying DDE \eqref{eq:ddef} affects only the boundary condition in $\mathcal{D}(\mathcal{A})$. Indeed, the latter represents the above mentioned rule of extension, while the differentiation action is the generator of (continuous) time translation. Accordingly, \eqref{eq:IVPDDEF} is equivalently recast as the \textit{abstract} IVP
\begin{equation}\label{eq:so}
    \left\{\setlength\arraycolsep{0.1em}
    \begin{array}{rcl}
    u'(t) &=& \mathcal{A}u(t), \quad t \geq 0,\\[1mm]
    u(0)  &=& \varphi,
    \end{array}\right.
\end{equation}
in the function space $\mathcal{X}$, where the equivalence relies on $u(t) = x_t$ (to be precise the equivalence holds when $\varphi\in\mathcal{D}(\mathcal{A})$, otherwise one can talk about \textit{mild} solutions).

Although the previous formulation \eqref{eq:ddef} is general, many practical applications involve DDEs with a finite number $k$ of constant discrete delays ${0<\tau_{1}<\cdots<\tau_{k}}$. Specifically, these equations are described as 
\begin{equation}\label{eq:ddediscrete}
x'(t)=f(x(t),x(t-\tau_{1}),\ldots,x(t-\tau_{k})),
\end{equation}
where, with slight abuse of notation, we now consider ${f:\mathbb{R}^{n(k+1)}\to\mathbb{R}^{n}}$. 
A concise overview of the specific DDEs of the form \eqref{eq:ddediscrete} that we use in Section~\ref{sec:2} as prototype models to test the introduced methodologies, is anticipated below. 
Note that each model is considered for a specific task in view of discussing the application of SINDy.

\begin{example}\label{ex:logistic}
The \textit{delay logistic equation}, also known as the Hutchinson equation \citep{Hutchinson1948}, is one of the most well-known DDEs. It is given by
\begin{equation}\label{eq:logistic}
x'(t) = r x(t) \left(1 - \frac{x(t - \tau)}{K} \right),
\end{equation}
where ${r > 0}$ is the intrinsic growth rate of the population, ${K > 0}$ is the carrying capacity and ${\tau > 0}$ represents the time delay in the feedback mechanism corresponding to internal competition.
The qualitative behavior of the solution of the corresponding IVP depends on $r$ and $\tau$: monotonic when $r\tau<1/\mathrm{e}$ and oscillatory when ${1/\mathrm{e}<r\tau<\pi/2}$. 
In both scenarios, the solution converges to the unique nontrivial equilibrium point ${\bar x=K}$ as ${t\to\infty}$ and $K$ also determines how fast the system converges. 
At the critical parameter value ${r\tau=\pi/2}$ a Hopf bifurcation \citep{kuznetsov1998elements} occurs and for ${r\tau>\pi/2}$ the solution tends to a stable limit cycle. 
The limit cycle can be represented using the Lambert $W$ function and the structure of the cycle varies as $r$ increases, resulting in sections with steep derivatives that are interspersed with plateaus~\citep{Gopalsamy1992}. 
In the sequel we use this model as a fundamental scalar example featuring a single discrete delay, making it an ideal starting point for demonstrating the effectiveness of sparse identification techniques for DDEs in the case the delay is assumed to be known.
\end{example}

\begin{example}\label{ex:mg}
The \textit{Mackey-Glass equation}~\citep{mg77} is the nonlinear DDE
\begin{equation}\label{eq:mg}
x'(t)=\beta\frac{x(t-\tau)}{1+(x(t-\tau))^{\alpha}}-\gamma x(t),
\end{equation}
with $\beta$, $\gamma$ and $\alpha$ positive real numbers. 
It was initially introduced to model the generation of blood cells in physiology, based on a feedback mechanism that relies on a delayed Hill-type response represented by the term
\begin{equation}\label{eq:hill}
h(t-\tau;\alpha)\coloneqq\frac{1}{1+(x(t-\tau))^{\alpha}},
\end{equation}
where $\alpha$ is the Hill coefficient. 
Depending on the parameter values, \eqref{eq:mg} can exhibit periodic solutions or possibly chaotic dynamics. 
As anticipated, \eqref{eq:mg} is particularly relevant due to its nonlinear rational term~\eqref{eq:hill}, which is representative of characteristic feedback in biological systems. 
Therefore, it will be explicitly incorporated into the framework for sparse identification to accurately extract the dynamics from data. 
In the sequel we use this model to optimize both $\tau$ and $\alpha$ assuming their values are unknown.  
\end{example}

\begin{example}\label{ex:multiple}
The \textit{two-neuron model} with multiple time delays
\begin{equation}\label{eq:neuron}
\left\{\setlength\arraycolsep{0.1em}\begin{array}{rcl}
x^\prime_1(t) &=& -\kappa x_1(t) + \beta \tanh(x_1(t-\tau_\mathrm{s})) + a_{12} \tanh(x_2(t-\tau_2)), 
\\[1mm]
x^\prime_2(t) &=& -\kappa x_2(t) + \beta \tanh(x_2(t-\tau_\mathrm{s})) + a_{21} \tanh(x_1(t-\tau_1)),
\end{array}\right.
\end{equation}
was originally introduced by \cite{shayer2000} to describe a simple neural network consisting of two interconnected neurons. Here, $x_1(t)$ and $x_2(t)$ represent the activity of the two neurons at time $t$, while ${\kappa > 0}$ denotes the decay rate, modeling the natural dissipation of neural activity over time. The parameter $\beta$ represents the self-feedback strength of each neuron through a delayed activation function and $a_{12}$, $a_{21}$ are the coupling coefficients, capturing the interaction strength between the neurons. 
The function $\tanh$ serves as an activation function that models the nonlinear response of the neurons. 
Finally, $\tau_\mathrm{s}$ represents the self-feedback delay, while $\tau_1,\tau_2$ denote the interaction delays between neurons. 
We consider this model as an advanced test case because it exemplifies a system of physical dimension ${n=2}$ with multiple discrete delays, thus providing a robust benchmark to assess the performance of sparse identification techniques in more complex scenarios. 
\end{example}

\CCLsubsection{Introducing SINDy for Ordinary Differential Equations}\label{sec:sindy}

SINDy~\citep{bpk16, bk19, champ19} is a data-driven framework for discovering governing equations of dynamical systems. 
By considering that most dynamical systems have an underlying sparse structure in terms of RHS description, SINDy constructs interpretable models from time-series data. 
This is particularly useful when the true governing equations are unknown or when data-driven model reduction is desired. 
The method assumes that the system state $x(t)$ is measured at discrete, either equidistant or random, time points ${t_1, \ldots, t_m}$ and its evolution is governed by the ODE
\begin{equation}\label{eq:ode}
    x' = f(x),
\end{equation} 
for an unknown RHS ${f:\mathbb{R}^{n}\to\mathbb{R}^{n}}$. 
To approximate the latter, SINDy expresses each component $f_i(x)$, ${i = 1, \ldots, n}$, as a sparse linear combination of $p$ candidate basis functions $\theta_j(x)$:
\begin{equation}\label{eq:fi}
f_i(x) = \sum_{j=1}^{p} \theta_j(x) \xi_{j,i}.
\end{equation}
Here, the functions ${\theta_1, \ldots, \theta_p}$ serve as building blocks for approximating the underlying system behavior and they are chosen in accordance with the expected characteristics of the system. 
They typically include polynomials, with possible additional terms such as \text{trigonometric functions}, \text{exponentials}, \text{rational} or other nonlinear functions. 
\smallskip

To represent the system reconstruction and estimate the coefficients ${\xi_{j,i} , j = 1, \ldots, p}$, for each ${i = 1, \ldots, n}$, we introduce the following notation. The measurement data are arranged in a matrix
\begin{equation*}
    X \coloneqq 
    \begin{bmatrix}
        x_1(t_1) &  \dots & x_n(t_1) \\      
        \vdots & \ddots  & \vdots \\ 
        x_1(t_m)  & \dots & x_n(t_m)
    \end{bmatrix}\in \mathbb{R}^{m \times n},
\end{equation*}
where rows correspond to time points and columns to state variables. Similarly, a matrix
\begin{equation*}
    X' \coloneqq 
    \begin{bmatrix} 
    x_1'(t_1)  & \dots & x_n'(t_1) \\
    \vdots & \ddots  & \vdots \\
    x_1'(t_m)  & \dots & x_n'(t_m)
    \end{bmatrix}\in \mathbb{R}^{m \times n},
\end{equation*}
is defined to collect the derivative values, which may be directly available from measurements or estimated using robust numerical differentiation techniques~\citep{chtd11}. 
To approximate the function $f$ in \eqref{eq:ode} according to \eqref{eq:fi}, a library matrix ${\Theta(X) \in \mathbb{R}^{m \times p}}$ is introduced, containing the $p$ candidate basis functions evaluated at all $m$ time points. 
As an example consider a polynomial library with additional trigonometric terms, e.g.,
\begin{equation*}
\Theta(X) \coloneqq \begin{bmatrix}
1 & X & X^2 & \cdots & X^d & \cdots & \sin(X) & \cos(X) & \cdots
\end{bmatrix}.
\end{equation*}
Above, $X^d$ denotes a matrix containing all possible polynomial combinations of degree $d$ in the components of $x$. For example, for $d = 2$,
\begin{equation*}
X^2 \coloneqq 
\begin{bmatrix}
x_1^2(t_1) & x_1(t_1)x_2(t_1) & \cdots & x_n^2(t_1) \\
x_1^2(t_2) & x_1(t_2)x_2(t_2) & \cdots & x_n^2(t_2) \\
\vdots & \vdots & \ddots & \vdots \\
x_1^2(t_m) & x_1(t_m)x_2(t_m) & \cdots & x_n^2(t_m)
\end{bmatrix}.
\end{equation*}
Similarly, for trigonometric functions,
\begin{equation*}
\sin(X) \coloneqq 
\begin{bmatrix}
\sin(x_1(t_1))  & \cdots & \sin(x_n(t_1)) \\
\vdots & \ddots & \vdots \\
\sin(x_1(t_m)) & \cdots & \sin(x_n(t_m))
\end{bmatrix}.
\end{equation*}
The choice of basis functions significantly impacts the success of the method. 
As discussed in~\citep{bpk16, bk19}, polynomial basis functions are commonly used for modeling ODEs, but other function types may be more suitable depending on the system. 
The number and type of basis functions determine the number $p$ of columns in $\Theta(X)$. 
To give a measure of the relevant complexity, the total number of basis functions of a polynomial library of degree $d$ for a system of physical dimension $n$ is ${p = \binom{n + d}{d}}$.

To determine the coefficients $\xi_{j,i}$ in \eqref{eq:fi} we assemble them into vectors ${\xi_i \in \mathbb{R}^p}$ for each ${i = 1, \ldots, n}$ and the overall problem is formulated as ${X' = \Theta(X)\Xi}$, where ${\Xi \in \mathbb{R}^{p \times n}}$ is the matrix ${\Xi=[\,\xi_1\; \xi_2\; \cdots\; \xi_n\, ]}$. 
The coefficients for each column $\xi_i$ are identified by solving independently the optimization problem
\begin{equation}\label{eq:optknown}
\xi_i =  \argmin_{\xi \in \mathbb{R}^p}\left(\left\|X_{i}'-\Theta(X)\xi\right\|_2+\lambda\left\|\xi\right\|_1\right) ,
\end{equation}
where $X_i'$ is the $i$-th column of $X'$. 
In the RHS of \eqref{eq:optknown} the first term (in $\|\cdot\|_{2}$) measures the reconstruction error while the second term (in $\|\cdot\|_{1}$) is a regularization promoting sparsity through the hyperparameter ${\lambda \geq  0}$. 
Methods such as sequential thresholded least squares (STLS)~\citep{bpk16} or least absolute shrinkage and selection operator (LASSO) regression~\citep{lasso96} can be employed to solve \eqref{eq:optknown}. 
The latter is an $\ell_{1}$ regression analysis method that performs both variable selection and regularization in order to enhance the prediction accuracy and interpretability of the produced model. 
The LASSO method minimizes the usual sum of squared residuals, with a constraint on the sum of the absolute values of the coefficients.
This constraint tends to produce some coefficients that are exactly zero, effectively performing variable selection. 
Alternatively, STLS is a method for tackling optimization problems that necessitate sparsity. 
In this approach, one repeatedly solves the least squares problem and then applies a thresholding operation: any coefficient whose absolute value is below a designated threshold (given by $\lambda$) is set to zero. 
This iterative process is carried out for a predetermined number of iterations, ultimately yielding a notably sparse solution. 
STLS is particularly effective in high-dimensional contexts where the number of features \(p\) significantly exceeds the number of observations \(m\)~\citep{bpk16}.

\CCLsubsection{{E-SINDy}: Expert SINDy for Known Delays}\label{sec:Etauknown}

Extending SINDy to DDEs~\citep{bbt24,kopeczi23,pec24,sandoz23} introduces additional complexities due to the system's dependence on past values. 
As discussed in Section~\ref{sec:DDEs}, we consider the case of DDEs with a finite number of constant discrete delays, specifically those of the form given by~\eqref{eq:ddediscrete}. 
When both the number and the values of the delays are known, the SINDy framework can be extended to handle such systems by incorporating delayed copies of the state variables into the library functions. 
Explicitly, the problem can be reformulated as 
\begin{equation}\label{eq:sindyaugmented}
X' = \Theta(X, X_{\tau})\Xi(\tau),
\end{equation}
where $X_{\tau}$ is a shorthand for $X_{\tau_{1}},\ldots,X_{\tau_{k}}$, with 
\begin{equation*}
    X_{\tau_{i}} \coloneqq 
    \begin{bmatrix}
        x_1(t_1-\tau_{i}) & \dots & x_n(t_1-\tau_{i}) \\
        \vdots & \ddots  & \vdots \\
        x_1(t_m-\tau_{i})  & \dots & x_n(t_m-\tau_{i})
    \end{bmatrix},
\end{equation*}
for $i=1,\ldots,k$ containing the delayed samples ${x(t-\tau_i)}$.
Note that including delays increases the dimensionality of a polynomial library, with the number of terms growing from ${p = \binom{n+d}{d}}$ to ${q = \binom{(k+1)n + d}{d}}$ when $k$ delays are considered with polynomials of degree $d$.

In earlier papers~\citep{sandoz23,pec24} it was assumed that the delayed samples can be directly retrieved from the measurements, based on the idea that each delay corresponds to an integer multiple of a uniform sampling period. 
Here we discard this assumption by acquiring delayed samples via piecewise linear interpolation of the available measurements. 
This allows for considering both incommensurate delays and non-uniformly distributed sampling times, thus enhancing the approach's applicability. 
The augmented library matrix ${\Theta(X, X_{\tau}) \in \mathbb{R}^{m \times q}}$ is then constructed by combining basis functions of both the current data and the delayed data, e.g.,
\begin{equation*}
\Theta(X, X_{\tau}) \coloneqq 
\begin{bmatrix}
1 & X & X_{\tau} & X^2 & XX_{\tau} & X_{\tau}^2 & \cdots & \sin(X) & \sin(X_{\tau}) & \cdots
\end{bmatrix}.
\end{equation*}

The sparse regression problem \eqref{eq:sindyaugmented} now depends on the known delay values ${\tau_{1},\ldots,\tau_{k}}$ and is solved again column-wise:
\begin{equation*}
\xi_i(\tau) = \argmin_{\xi \in \mathbb{R}^q}\left(\|X_i' - \Theta(X, X_{\tau})\xi\|_2 + \lambda \|\xi\|_1\right).
\end{equation*}
The notation $\Xi(\tau)$ in \eqref{eq:sindyaugmented} highlights the dependence of the sparse solution $\Xi$ on the given vector of delay values ${\tau:=[\,\tau_{1},\ldots,\tau_{k}\,]^{\top}}$.
For later use the corresponding reconstruction error is denoted as
\begin{equation}\label{epsilontau}
\epsilon(\tau) := \|X' - \Theta(X, X_{\tau})\Xi(\tau)\|_{2}.
\end{equation}

As discussed earlier, this extension of SINDy to DDEs of the form \eqref{eq:ddediscrete} assumes that both the number and values of the delays are known a priori, which is uncommon in practical applications and often relies on prior modeling expertise. 
The latter observation justifies the name \textit{{E-SINDy}} in \citep{bbt24}, where ``E'' stands for ``Expert''\footnote{\label{note2} The term ``expert'' reflects an approach that directly incorporates prior knowledge of delays to handle DDEs, and in general it tackles a DDE as it is. 
This terminology follows the convention in \citep{bdls16}.}. 
Thus this methodology leverages prior knowledge of delays to directly recover the concerned DDE. 
Moreover, {E-SINDy} primarily aims to reconstruct the RHS $f$, which is essential for tasks such as model interpretation rather than long-term predictive simulations. 
However, it is important to emphasize that accurately reproducing the true trajectory of the system does not necessarily require the correct identification of the RHS function. 
Indeed, even if the reconstructed model successfully replicates the observed data, it may still fail to capture the true governing equations of the system \citep{pec24}. 

\begin{table}[htbp]
\centering
\small 
\begin{tabular}{lrr}
 \toprule
  & \multicolumn{2}{c}{Samples Type} \\
 \cmidrule(r){2-3}
 & Uniform & Random \\
 \midrule
 Samples Size ($m$) & 10 & 10 \\
 \midrule
 \multicolumn{3}{l}{\textbf{Training RMSE$_{x^\prime}$}} \\
 STLS & \boldmath{$7.80\times10^{-15}$} & \boldmath{$4.19\times10^{-15}$} \\
 LASSO & $2.80\times10^{-3}$ & $2.66\times10^{-3}$ \\
 \midrule
 \multicolumn{3}{l}{\textbf{Testing RMSE$_{x^\prime}$}} \\
 STLS & \boldmath{$4.62\times10^{-15}$} & \boldmath{$6.43\times10^{-15}$} \\
 LASSO & $2.61\times10^{-3}$ & $2.89\times10^{-3}$ \\
 \midrule
 \multicolumn{3}{l}{\textbf{Testing RMSE$_{x}$}} \\
 STLS & \boldmath{$4.93\times10^{-14}$} & \boldmath{$6.09\times10^{-14}$} \\
 LASSO & $6.35\times10^{-1}$ & $2.04\times10^{-2}$ \\
 \bottomrule
 \end{tabular}
\caption{Performance comparison of {E-SINDy} on \eqref{eq:logistic} using STLS and LASSO methods with uniform and random samples, showing RMSE values for training and testing phases; see text for more details.}
 \label{tab:comparisondl}
 \end{table}

To illustrate the application of {E-SINDy} to DDEs, we reconstruct the delay logistic equation \eqref{eq:logistic} for ${r=1.8}$, ${K=10}$ and ${\tau=1}$ assumed to be known. 
We use the time series obtained by solving the IVP for \eqref{eq:logistic} with history function ${\varphi(s) = \cos(s)}$ (via MATLAB's \texttt{dde23}) on the time window ${[0,30]}$, using the data in ${[0,18]}$ for training {E-SINDy} and the data in ${[18,30]}$ for testing the obtained model. 
A polynomial library of degree ${d=2}$ is employed. 
To evaluate the reliability of the method we consider both uniform and random (non-uniform) distributions for the ${m = 10}$ data points used for training. 
We assume that at these data points both the state $x$ as well as the derivative $x'$ are accurately measured.
Results are collected in Table~\ref{tab:comparisondl}, where the following root mean squared errors are considered:
\begin{equation}\label{RMSEx'E}
\textrm{RMSE}_{x'} := \|X' - \Theta(X, X_{\tau})\Xi(\tau)\|_{2} ,
\end{equation}
\begin{equation}\label{RMSExE}
\textrm{RMSE}_{x} :=\|X - \mathcal{S}_{\text{DDE}}(\Theta(X, X_{\tau})\Xi(\tau))\|_{2}.
\end{equation}
The former is indeed $\epsilon(\tau)$ defined in \eqref{epsilontau}, while in the latter the notation $\mathcal{S}_{\text{DDE}}$ represents the DDE solver that progresses the model forward in time (e.g., MATLAB's \texttt{dde23} as we adopt in this work). 
The results show the efficacy of {E-SINDy} in reconstructing the RHS (both training and testing), as well as how closely the reconstructed trajectory matches the true solution (only testing). 
Tests concerned both STLS and LASSO, showing the superiority of the former. 

\begin{figure}[!t]
\centering
\includegraphics[width=1\linewidth]{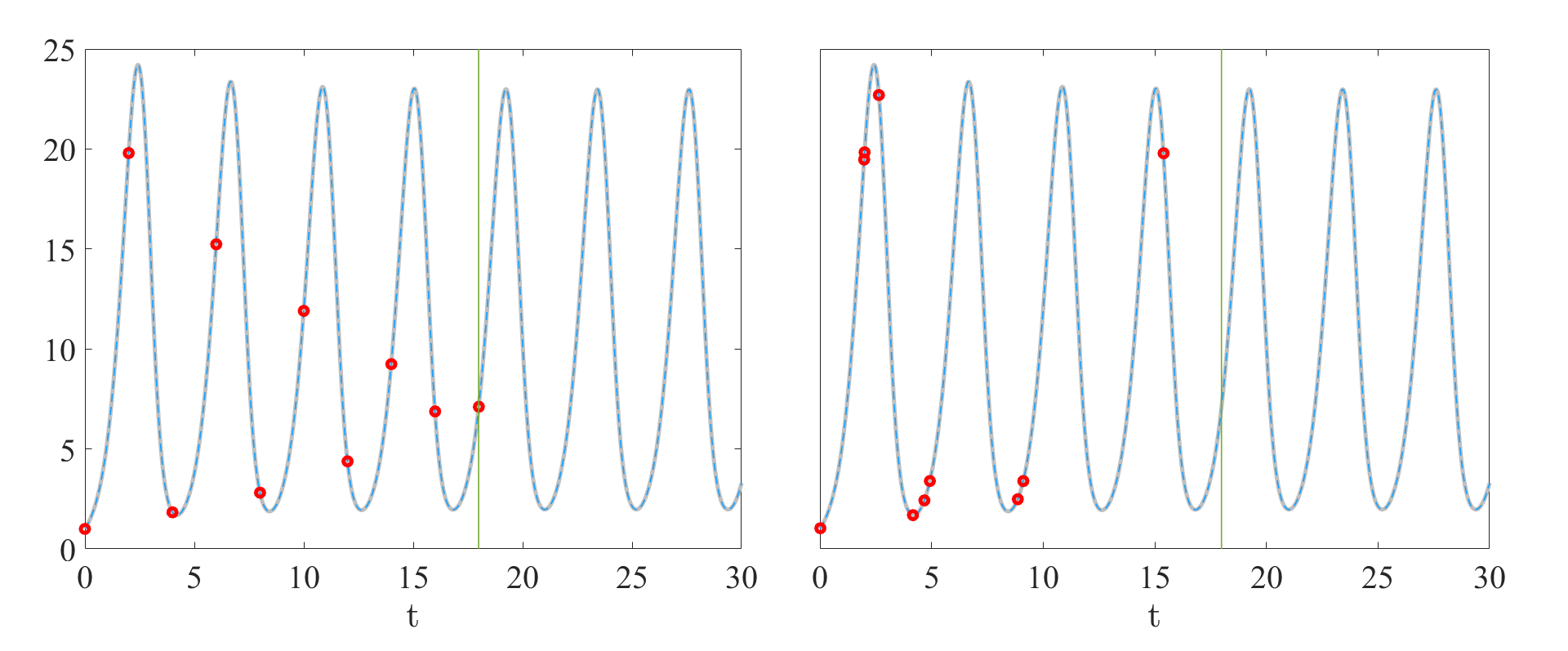}
\caption{True (gray) and {E-SINDy} trajectory (cyan) of the delay logistic equation \eqref{eq:logistic} for ${r=1.8}$, ${K=10}$ and ${\tau=1}$, reconstructed with ${m=10}$ uniform (left) and random (right) training samples (red dots) in ${[0,18]}$, and tested in ${[18,30]}$.}\label{fig:DLE_known} 
\end{figure}

Figure~\ref{fig:DLE_known} presents the reconstructed trajectories for both uniform (left) and random (right) training samples, demonstrating the effectiveness of {E-SINDy} in identifying the governing dynamics of \eqref{eq:logistic}.
Observe that in general very few samples are enough to correctly reconstruct the underlying model.

\CCLsubsection{Unknown Delays and Parameters}\label{sec:Etauunknown}

In practical scenarios, the number and values of the delay parameters in a DDE are often unknown. 
In the context of SINDy, one way to handle this challenge is to consider minimizing the map ${\tau \mapsto \epsilon(\tau)}$ of the reconstruction error \eqref{epsilontau} or \eqref{RMSEx'E}. 
When only a single delay is unknown \cite{sandoz23} suggests a simple brute force (BF) approach~\citep{brute2011,brute2019}: evaluate $\epsilon(\tau)$ for a range of candidate delay values and pick the one that minimizes the reconstruction error. 
Although straightforward, this approach can be computationally expensive if the search space is large: the number of calls to SINDy is the number of considered candidate delays. 

A more sophisticated method is suggested by \cite{pec24}, in which Bayesian optimization (BO)~\citep{williams2006,shahriari2015,garnett2023} is employed to identify the delay value that minimizes $\epsilon(\tau)$, decreasing the number of calls to {E-SINDy}.
In addition, \cite{pec24} extended this idea to a multivariate scenario, simultaneously optimizing multiple delays and/or potentially other unknown non-multiplicative parameters (e.g., non-multiplicative parameter entries in the function library). 
This extension is especially beneficial for cases like the Mackey-Glass equation \eqref{eq:mg}, where the presence of the nonlinear rational term \eqref{eq:hill} requires identifying both $\tau$ and $\alpha$.

Building on these methods, \cite{bbt24} presented the alternative optimization based on particle swarm (PS)~\citep{bonyadi2017,kennedy1995,shi1998}, which explores a continuous range for the unknown delay(s) or parameter(s), offering in general a more efficient alternative to BO and BF, particularly in higher-dimensional parameter spaces.

To summarize, let us stress that ``external optimization'' methods identify unknown delays and parameters, while ``internal optimization'' in SINDy (e.g., STLS \citep{bpk16}) performs sparse regression to determine the active terms in the model based on these estimated parameters.

\begin{figure}[!t]
\centering
\includegraphics[width=1\linewidth]{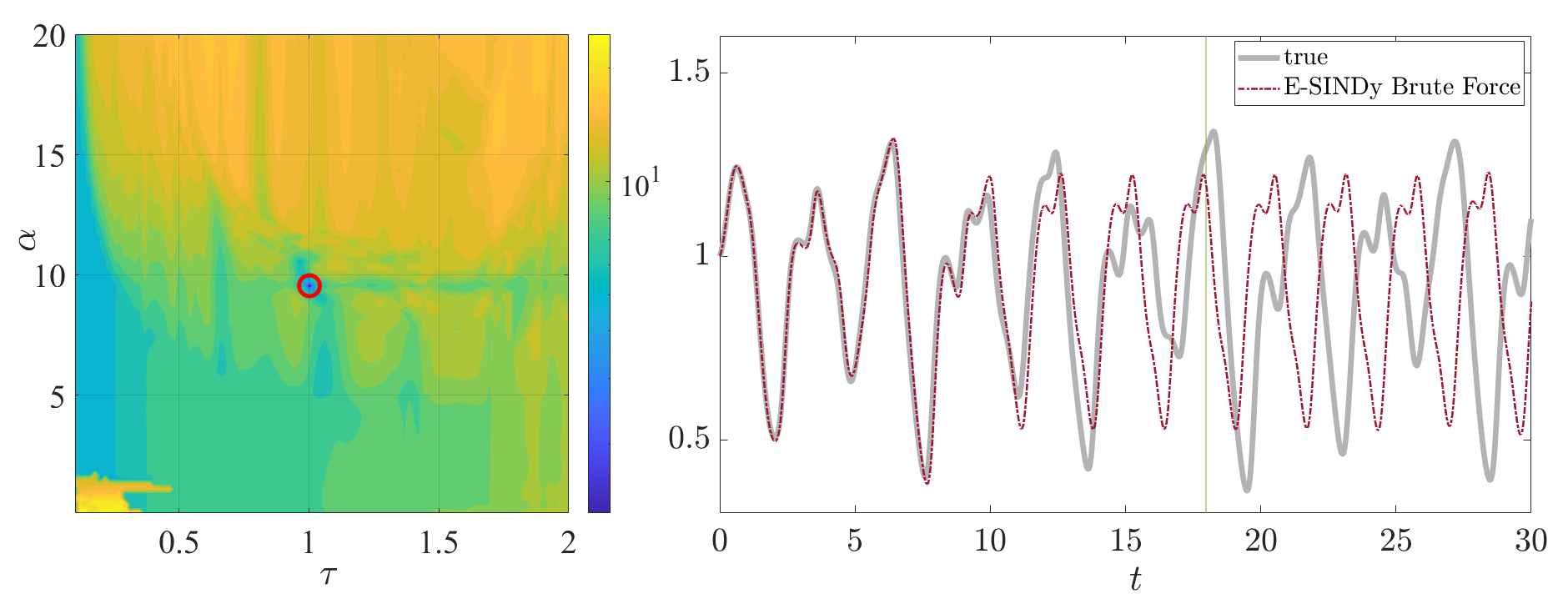}
\caption{(left) External optimization of the unknown delay $\tau$ and the Hill exponent $\alpha$ for the Mackey-Glass equation \eqref{eq:mg} using BF.
Color indicates the ${\textrm{RMSE}_{x'}}$ over a ${100 \times 100}$ grid in the ${(\tau,\alpha)}$ plane, the nominal values are ${(\tau,\alpha)=(1, 9.6)}$, and the optimal point ${(\widehat{\tau},\widehat{\alpha})=(1.0020, 9.5475)}$ is marked by red circle. 
(right) True (gray) and {E-SINDy} (red) trajectories after external optimization.}\label{fig:MG_BF}
\end{figure}

We conclude this section by testing {E-SINDy} with external optimization to recover both $\tau$ and $\alpha$ for the Mackey-Glass equation~\eqref{eq:mg}. 
Data are obtained in the time window ${[0,30]}$ by simulating the corresponding IVP with initial function ${\varphi(s) = \cos(s)}$ and reference parameters ${\beta = 4}$, ${\gamma = 2}$ and nominal values ${\tau = 1}$ and ${\alpha = 9.6}$ for the unknown delay and Hill exponent. 
We use ${m = 100}$ training samples in ${[0,18]}$ and the remaining interval ${[18,30]}$ for testing. 
The chosen polynomial library has degree ${d=2}$ and includes the rational term \eqref{eq:hill} with non-multiplicative unknown parameter $\alpha$ affecting the model's nonlinearity.

Figure~\ref{fig:MG_BF} (left) shows the reconstruction error \eqref{RMSEx'E} as a function of both $\tau$ and $\alpha$ using a ${100 \times 100}$ grid in the interval ${[0.1, 2]}$ for $\tau$ and ${[0.1, 20]}$ for $\alpha$.
Figure~\ref{fig:MG_BF} (right) shows the corresponding trajectory obtained with the model reconstructed by {E-SINDy} externally optimizing $\tau$ and $\alpha$ via BF. 
The obtained values are ${[\widehat{\tau},\widehat{\alpha}]=[1.0020, 9.5475]}$. 
The discrepancy of the reconstructed trajectory is due to the low accuracy of the final optimal values. 
Better results are summarized in Table \ref{tab:comparison2} of Section \ref{sec:discussion} for BO and PS following a more detailed comparison.

A significant drawback of {E-SINDy} is represented by the assumption that the number of delays is assumed to be known, or at least the adopted library should allow for a sufficiently large number of delayed terms. 
As such information is often unavailable in practice, the methodology illustrated in the forthcoming section aims at addressing the issue,  operating without requiring prior knowledge of the number of delays.

\CCLsubsection{{P-SINDy}: Pragmatic SINDy via Pseudospectral Collocation}\label{sec:PSINDy}

{E-SINDy} for unknown delays and parameters as introduced in Section~ \ref{sec:Etauunknown} relies on externally optimizing (via BF, BO or PS) a predefined number $\bar{k}$ of delays, assuming ${\bar{k} \geq k}$, where $k$ is the true number of delays in \eqref{eq:ddediscrete}. 
If ${\bar{k} > k}$ and {E-SINDy} succeeds, the resulting sparse matrix $\Xi$ will contain ${\bar{k}-k}$ null coefficients. 
However, multivariate optimization can be computationally expensive. 
To address this issue \cite{bbt24} introduced a more efficient approach where only the maximum delay in \eqref{eq:ddediscrete} is required to be externally optimized. 
This method, termed \textit{{P-SINDy}} where P stands for ``Pragmatic''\footnote{\label{note3}This terminology follows again the naming convention in \citep{bdls16} where ``pragmatic'' refers to first approximating a DDE with an ODE and then applying techniques for ODEs.}, 
focuses first on reducing the problem to a finite-dimensional system of ODEs through pseudospectral collocation \citep{bdgsv16}, and then applies SINDy to the resulting ODE. 
Unlike {E-SINDy}, which relies on prior knowledge of delay parameters to reconstruct an interpretable RHS, {P-SINDy} adopts a ``black-box'' approach, aiming only at matching the observed trajectory and then simulating beyond the measurements time span. 
While {P-SINDy} sacrifices the interpretability of the reconstructed RHS, it requires no prior information on the number or values of delays and thus achieves enhanced computational efficiency by possibly optimizing only the maximum delay ${\bar{\tau}\coloneqq\tau_{k}}$ in the case the latter were unknown. 

Note that in {P-SINDy} the infinite-dimensional dynamical system \eqref{eq:so} on the state space $\mathcal{X}$ is reduced to a system of ODEs, each of which is itself a system of dimension $n$, corresponding to the original physical dimension. 
In the meantime, the actual dimension of the sparse regression problem remains confined to  $n$, similarly to {E-SINDy}, as only the first block of ODEs needs to be reconstructed, as illustrated below. 
Now we recall the necessary ingredients of pseudospectral collocation and then describe how to augment the SINDy library to obtain {P-SINDy}.

\CCLsubsubsection{Pseudospectral Collocation}\label{sec:PSC}

Consider the IVP \eqref{eq:IVPDDEF} and recall \eqref{eq:generator} and \eqref{eq:so}. 
To approximate \eqref{eq:IVPDDEF} with a finite-dimensional system of ODEs, we apply a \emph{pseudospectral collocation} scheme based on \emph{Chebyshev polynomials} \citep{breda2022controlling}. 
This method provides a systematic way to reduce the infinite-dimensional state space $\mathcal{X}$ to a finite-dimensional subspace $\mathcal{X}_M$ and, correspondingly, a DDE to a system of ODEs, which can then be treated with relevant tools (as the standard SINDy for ODEs).

For a positive integer $M$ let
\begin{equation*}
s_{i}:=\frac{\bar\tau}{2}\left(\cos\left(\frac{i\pi}{M}\right)-1 \right),\quad i=0,1,\ldots,M,
\end{equation*}
be the ${M+1}$ \emph{Chebyshev extremal nodes} over ${[-\bar{\tau},0]}$. 
Define ${\mathcal{X}_M :=\mathbb{R}^{n(M+1)}}$ as the discrete counterpart of $\mathcal{X}$ and introduce the following \emph{restriction} and \emph{prolongation} operators:
\begin{equation*}
\setlength\arraycolsep{0.1em}\begin{array}{ll}
R_{M}:\mathcal{X}\to\mathcal{X}_{M},&\quad R_{M}\psi:=[\,\psi(s_0),\psi(s_1),\ldots,\psi(s_M)\,]^{\top},
\\[1mm]
P_{M}:\mathcal{X}_{M}\to\mathcal{X},&\quad P_{M}\Psi:=\sum_{j=0}^{M}\ell_j(s)\Psi_j.
\end{array}
\end{equation*}
The former evaluates a function ${\psi \in \mathcal{X}}$ at the collocation nodes, while the latter reconstructs a continuous function from its discrete values via Lagrange interpolation. The Lagrange basis functions ${\{\ell_j\}_{j=0}^{M}}$ are defined on the collocation nodes ${\{s_j\}_{j=0}^{M}}$ as
\begin{equation*}
    \ell_j(s) := \prod_{\substack{k=0 \\ k\neq j}}^{M} \frac{s - s_k}{s_j - s_k}, \quad j=0,1,\dots,M, \quad s \in[-\bar{\tau},0].
\end{equation*}

In the discrete setting, the infinitesimal generator $\mathcal{A}$ is approximated by the finite-dimensional operator ${\mathcal{A}_{M}:\mathcal{X}_{M}\to\mathcal{X}_{M}}$ given by
\begin{equation*}
[\mathcal{A}_M\Psi]_0=f(P_M\Psi),\qquad[\mathcal{A}_M\Psi]_i= \sum_{j=0}^{M} d_{i,j} \Psi_j,\quad i=1,\ldots,M,
\end{equation*}
where the derivatives ${d_{i,j}\coloneqq \ell_j'(s_i)\, I_n}$, ${i=1,\dots,M}$, ${j=0,1,\dots,M}$, of the Lagrange basis functions can be calculated explicitly and collected as entries of the so-called Chebyshev differentiation matrix \citep{trefethen2000}. 
The discrete expressions of $\mathcal{A}_{M}$ above correspond to, respectively, the domain condition and the action of $\mathcal{A}$ in \eqref{eq:generator}. 
Note that the RHS $f$ of the original DDE affects only the first ODE, which is indeed a block of $n$ ODEs. 
In this respect, $\mathcal{A}_{M}$ faithfully reproduces the essential features of $\mathcal{A}$.

The above discretization leads to the following finite-dimensional IVP approximating \eqref{eq:IVPDDEF}:
\begin{equation}\label{DDEFM}
\left\{\setlength\arraycolsep{0.1em}\begin{array}{rcll}
U_{0}'(t) &=& f\bigl(P_M\,U(t)\bigr), & \\[3pt]
U_{i}'(t) &=& D_M\,U(t), & \quad i = 1,\ldots,M,\\[3pt]
U_{i}(0) &=& \varphi\bigl(s_{i}\bigr), & \quad i = 0,1,\ldots,M,
\end{array}\right.
\end{equation}
where
\begin{equation*}
      D_M\coloneqq
  \begin{bmatrix}
  \ell_0'(s_1)  & \dots & \ell_M'(s_1) \\
  \vdots               & \ddots & \vdots          \\
  \ell_0'(s_M)  & \dots & \ell_M'(s_M)
  \end{bmatrix},
\end{equation*}
${U(t)\coloneqq[\,U_{0}(t),U_{1}(t),\ldots,U_{M}(t)\,]^{\top}\in\mathcal{X}_{M}}$, and each component $U_i(t)$ is to approximate ${x_t(s_i) = x\bigl(t+s_i\bigr)}$.  
As anticipated, ${U_0(t)\approx x(t)}$ is the only quantity of interest and the corresponding ODE is the only one reconstructed by SINDy. 
However, it is coupled to all the other components ${U_i(t), i=1,\ldots,M}$, through the interpolation polynomial $P_MU(t)$.

To summarize, the pseudospectral collocation approach, leveraging Chebyshev nodes and the associated differentiation matrix, typically achieves spectral accuracy \citep{trefethen2000,bary04}. 
For delay problems, it provides a systematic way to approximate infinite-dimensional dynamics within a finite-dimensional ODE framework. 
Note that alternative discretization approaches may also be implemented, such as those described in \citep{lehotzky2016pseudospectral,borgioli2020pseudospectral}.

\CCLsubsubsection{Extended Collocation Library}\label{sec:extlib}

As already observed, most of the ODEs in \eqref{DDEFM} (viz.\ the $n$-dimensional blocks for ${i=1, \ldots, M}$) arise from the differentiation action of $\mathcal{A}$, and are thus independent of the specific RHS $f$. 
The maximum delay $\bar\tau$ only affects the length of the collocation domain ${[-\bar\tau,0]}$. 
As the only ODE directly affected by $f$ is the first block ${U_0'(t) = f(P_M U(t))}$, {P-SINDy} applies sparse regression solely to this block. 
The resulting regression problem reads ${
X' = \Theta(X, X_s) \Xi_M(\bar{\tau})}$, where now $X_s$ denotes the collocated samples ${X_{s_1}, \ldots, X_{s_M}}$ with
\begin{equation*}
    X_{s_{i}} \coloneqq 
    \begin{bmatrix}
        x_1(t_1 + s_{i}) & \dots & x_n(t_1 + s_{i}) \\
        \vdots & \ddots & \vdots \\
        x_1(t_m + s_{i}) & \dots & x_n(t_m + s_{i})
    \end{bmatrix},
\end{equation*}
for ${i=1,\ldots,M}$ and the notation $\Xi_M(\bar{\tau})$ highlights the dependence of the sparse solution on the maximum delay $\bar{\tau}$ and on the collocation degree $M$.

As a result, {P-SINDy} operates without requiring the knowledge of possible intermediate delays.
Indeed, when ${k>1}$, the effect of multiple delays is mixed through interpolation, making difficult to recover the RHS $f$ in an interpretable form. 
Instead of directly identifying the RHS, {P-SINDy} prioritizes accurate trajectory reconstruction, reducing the problem to optimizing only the maximum delay $\bar{\tau}$, as opposed to $k$ delays in {E-SINDy}. 
This convenient pragmatic approach leverages well-established tools for ODEs, making it a computationally efficient alternative. 
Note finally that the size of a corresponding polynomial library of degree $d$ becomes ${r = \binom{(M+1)n + d}{d}}$.

We conclude this section with tests on the two-neuron system \eqref{eq:neuron} of 2 DDEs with 3 delays. 
For the reference parameters we utilize 
${\kappa = 0.5}$, ${\beta = -1}$, ${a_{12} = 1}$ and ${a_{21} = 2}$, 
with unknown delays having the nominal values 
${\tau_\mathrm{s} = 1.5}$, ${\tau_1 = 2}$ and ${\tau_2 = 2}$. 
The corresponding IVP with constant initial history ${\varphi(s) \equiv [\,0.5, -0.5 \,]^{\top}}$ is simulated on the time window ${[0,30]}$, producing ${m = 100}$ uniform samples in ${[0,18]}$ for training. 
We use a polynomial library of degree ${d=2}$ along with trigonometric functions: the polynomial terms address generic nonlinearities, while the trigonometric functions improve the modeling of oscillatory phenomena, as frequently seen in neural systems. 
In {E-SINDy} with BO the delay $\widehat{\tau}_s$ is optimized within the range ${[1, 2]}$ and the delays $\widehat{\tau}_1$ and $\widehat{\tau}_2$ are optimized within the range of ${[1, 2.5]}$. 
In {P-SINDy} with PS only ${\bar{\tau}=\widehat{\tau}_2}$ is optimized, using the range indicated above. 

\begin{figure}[!t]
\centering
\includegraphics[width=1\linewidth]{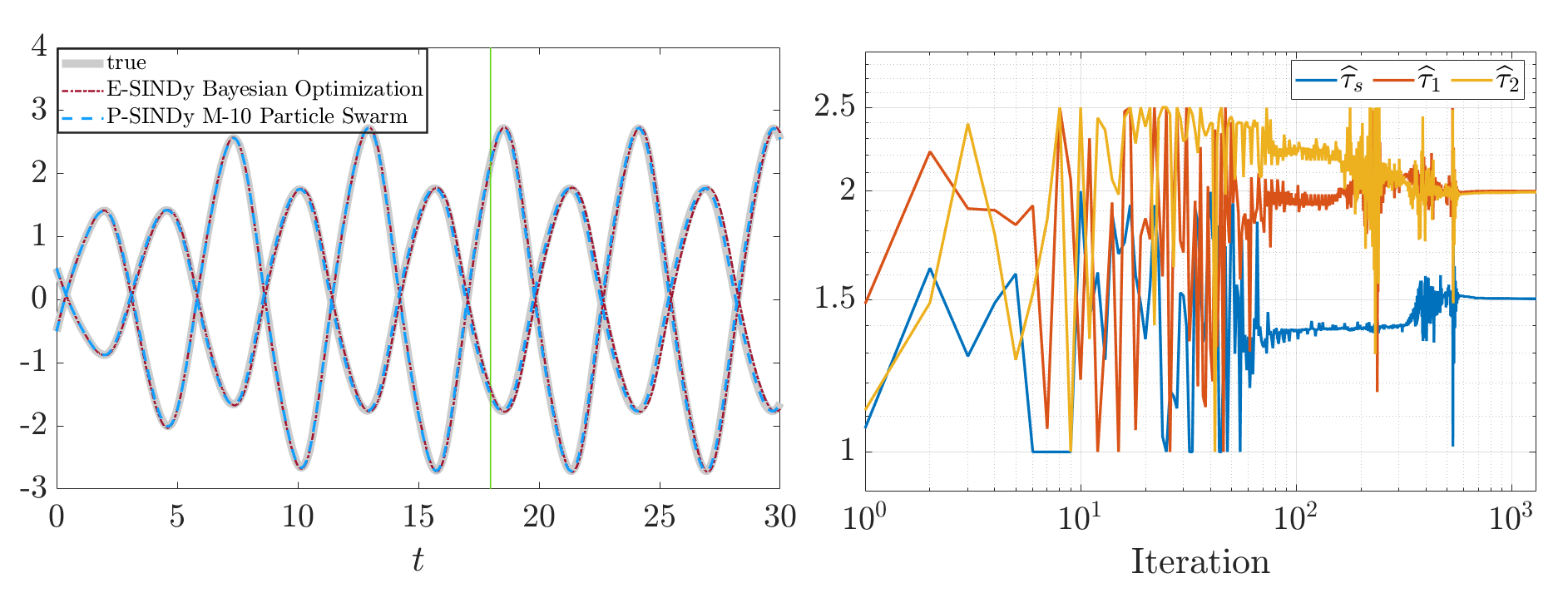}
\caption{(left) Reconstructed trajectory for the two-neuron model \eqref{eq:neuron} showing the true trajectory (gray), {E-SINDy} trajectory obtained with BO (red), {P-SINDy} trajectory of collocation degree ${M=10}$ obtained with PS (blue). 
The true delay values are ${(\tau_\mathrm{s},\tau_1,\tau_2)=(1.5, 2, 2)}$, the optimized delays values for {E-SINDy} are ${(\widehat{\tau}_s,\widehat{\tau}_1,\widehat{\tau}_2)=(1.5136, 1.9478,2.0598)}$ while for the {P-SINDy} approach the recovered delay is ${\bar{\tau}= 2.0019}$. 
(right) Delay recovery of delays for {E-SINDy} ${(\widehat{\tau}_{s},\widehat{\tau}_{1},\widehat{\tau}_{2})}$ using PS.}\label{fig:neuron} 
\end{figure}

Figure~\ref{fig:neuron} (left) illustrates the reconstructed trajectories for {E-SINDy} with BO and {P-SINDy} with PS and ${M=10}$. 
The results for {E-SINDy} with PS are omitted to avoid redundancy, as the trajectories are very similar to those obtained with {P-SINDy}. 
Instead, Figure~\ref{fig:neuron} (right) illustrates the evaluation cycle of the unknown delays obtained using PS with {E-SINDy}. 
The process is plotted against the number of iterations, basically representing the number of times {E-SINDy} is called. 
Further results are given in Table~\ref{tab:comparison1} in the forthcoming section. 
Additional test results with {P-SINDy} are presented in Section~\ref{sec:compare} where the following root mean squared errors are considered:
\begin{equation}\label{RMSEx'P}
\textrm{RMSE}_{x'} := \|X' - \Theta(X, X_{s})\Xi_{M}(\bar\tau)\|_{2}, 
\end{equation}
\begin{equation}\label{RMSExP}
\textrm{RMSE}_{x} :=\|X - \mathcal{S}_{\text{ODE}}(\Theta(X, X_{s})\Xi_{M}(\bar\tau))\|_{2}.
\end{equation}
Above, $\mathcal{S}_{\text{ODE}}$ computes the trajectory of the full system \eqref{DDEFM}, integrating the reconstructed first block of $n$ ODEs with the addition of the known differentiation terms from \eqref{DDEFM} (in this work we use MATLAB's \texttt{ode45}).

\CCLsubsection{Further Tests and Discussion on {E-SINDy} and {P-SINDy}}\label{sec:discussion}

Scope of this section is to provide further experiments on the DDEs \eqref{eq:logistic}, \eqref{eq:mg} and \eqref{eq:neuron} to evaluate the performance of both {E-SINDy} and {P-SINDy}, additionally employing BF, BO and PS as external optimizers for unknown delays or parameters. 
All the  tests were run on a Windows 11 system (CPU 5 GHz, RAM 16 GB) using MATLAB implementations for SINDy (available at \url{http://cdlab.uniud.it/software}, MATLAB version R2024a) and MATLAB's integrated optimizers, viz.\ \texttt{bayesopt.m} for BO from the Statistics and Machine Learning Toolbox and \texttt{particleswarm.m} for PS from the Global Optimization Toolbox for PS\footnote{To run {E-SINDy} open \texttt{ESINDy\_main.m} and modify the ``User inputs'' section according to desired model, time window, sample distribution and noise amplitude and then execute \texttt{ESINDy\_main} from the Command Window. 
Similarly, for {P-SINDy} open \texttt{PSINDy\_main.m}, adjust the user parameters (such as the collocation degree $M$) in the ``User inputs'' section and run \texttt{PSINDy\_main}. 
Each script automatically generates data for the selected model, constructs the library of candidate functions, optionally applies the chosen optimizer to estimate unknown delays and parameters and finally carries out the sparse regression step (via STLS).}. 

Data for each experiment consist of $m$ uniformly spaced samples obtained by solving the corresponding DDE IVP~\eqref{eq:IVPDDEF} via MATLAB’s \texttt{dde23} over the interval ${[0,T]}$ with a given history function $\varphi$. 
In each case 60\% of the samples are used for training and 40\% for testing. 
In the following tables and figures we label each SINDy approach as ``E'' for {E-SINDy} and ``P'' for {P-SINDy}, specifying the collocation degree as P5, P10 and P15 for respectively ${M=5,\,10}$ and $15$. 
Moreover, recovered values for DDE coefficients as well as for possible unknown delays and parameters are denoted with a hat.
The tables report the relevant absolute errors rather than the values themselves (e.g., ${|\tau-\hat\tau|}$ rather than $\hat\tau$), to better highlight the accuracy of the various reconstructions. 
As for the reconstruction errors, we refer to \eqref{RMSEx'E} and  \eqref{RMSExE} for {E-SINDy} and to \eqref{RMSEx'P} and  \eqref{RMSExP} for {P-SINDy}. 
Finally, note that the use of uniformly spaced samples do not restrict the values of unknown delays, as delayed and collocation samples are reconstructed via piecewise linear interpolation if not available.

A first series of tests is run on \eqref{eq:logistic} assuming by default ${K=1}$. The growth parameter and delay nominal values are set to ${r=1.8}$ and ${\tau=1}$, the latter assumed to be unknown and thus to be externally optimized. 
Furthermore, we use ${m=100}$ and ${T=30}$, with initial function ${\varphi(s)=\cos(s)}$, ${s\in[-1,0]}$, and a polynomial library with degree ${d=2}$. 
The latter describes the RHS through the vector ${\Xi=[\, 0,r,0,0,-r,0\,]^{\top}}$ corresponding to the library of candidate functions
\begin{equation}
\Theta(X,X_{\tau})=
\begin{bmatrix}
1 & X & X_{\tau} & X^{2} & XX_{\tau} & X^{2}_{\tau}
\end{bmatrix}.
\end{equation}
Table~\ref{tab:logistic} summarizes the results for E, P5, P10 and P15 with BF, BO and PS, including the absolute errors on both $r$ (as a measure of the goodness of the recovered sparse vector $\Xi$) and $\tau$, as well as the reconstruction errors $\textrm{RMSE}_{x'}$ and $\textrm{RMSE}_{x}$. 
In all the cases the external optimizers BF, BO and PS successfully identify the correct delay $\tau$ yet with slightly varying accuracy, with the precision of BF depending on the density of its candidate set (here $1\,000$ values in the range ${[0.1 , 2]}$). 
Note moreover that BO halts after a predetermined number of calls to SINDy (here 300), while PS terminates once a desired accuracy threshold is reached (here set to $10^{-3}$). 
As both {E-SINDy} and {P-SINDy} with different optimizers recover the coefficients and delay within good accuracy, the resulting trajectories for training and testing are nearly identical and we omit to plot them as the reader can refer back to Figure~\ref{fig:DLE_known} in Section \ref{sec:Etauknown}. 
The number of calls to SINDy and corresponding CPU times are reported in the first block of rows of Table~\ref{tab:comparison1}. 
Note how PS requires in general less calls to SINDy and consequently a much lower CPU time yet giving the same or even better accuracy, and that the number of calls to SINDy decreases for {P-SINDy} with increasing $M$, giving in general better accuracy yet higher CPU time. 
As compared to {E-SINDy} and BF, P5 or P10 with PS appear to be the most effective choices.

\begin{sidewaystable}[htbp]
\centering
\begin{tabular}{llrrrr}
\toprule
SINDy&optimizer&$|r-\widehat{r}|$&$|\tau-\widehat{\tau}|$ &RMSE$_{x^\prime}$& RMSE$_{x}$\\
\midrule
E&BF&\boldmath{$1.00\times10^{-5}$}& \boldmath{$3.00\times10^{-4}$}&  $1.72\times10^{-12}$ & $6.02\times10^{-5}$
\\
P5&BF& $1.00\times10^{-4}$&\boldmath{$3.00\times10^{-4}$}& $6.38\times10^{-10}$ & $2.43\times10^{-5}$\\
P10&BF&\boldmath{ $1.00\times10^{-5}$}&\boldmath{$3.00\times10^{-4}$}& $1.98\times10^{-13}$ &\boldmath{$1.25\times10^{-6}$}\\
P15&BF&\boldmath{$1.00\times10^{-5}$}&\boldmath{$3.00\times10^{-4}$}& \boldmath{$1.10\times10^{-14}$} & $8.83\times10^{-6}$\\
\midrule
E&BO& $1.00\times10^{-4}$ & $2.00\times10^{-4}$&  $3.72\times10^{-11}$  & $5.62\times10^{-3}$\\
P5&BO& $1.00\times10^{-4}$ & $1.00\times10^{-4}$ &  $2.68\times10^{-9}$  & $4.45\times10^{-3}$\\
P10&BO&\boldmath{$1.00\times10^{-5}$} & $1.00\times10^{-5}$ & $9.44\times10^{-11}$ & $3.48\times10^{-6}$\\
P15&BO&\boldmath{$1.00\times10^{-5}$}&\boldmath{$1.00\times10^{-6}$}&  \boldmath{$5.32\times10^{-12}$} &\boldmath{$6.35\times10^{-7}$}\\
\midrule
E&PS& $1.00\times10^{-4}$ &\boldmath{$1.00\times10^{-6}$} &  $4.91\times10^{-12}$  & $8.17\times10^{-5}$\\
P5&PS& $1.00\times10^{-2}$ & $1.00\times10^{-4}$ &   $1.72\times10^{-11}$  & $2.90\times10^{-5}$ \\
P10&PS& \boldmath{$1.00\times10^{-6}$} & $1.00\times10^{-5}$&   $2.38\times10^{-13}$ & $5.14\times10^{-7}$\\
P15&PS&\boldmath{$1.00\times10^{-6}$}& $1.00\times10^{-5}$ &  \boldmath{$4.01\times10^{-16}$} &\boldmath{$4.01\times10^{-8}$}\\
\bottomrule
\end{tabular}%
\caption{Performance comparison of {E-SINDy} (E) and {P-SINDy} (P, with ${M=5}$, ${10}$ and ${15}$) for \eqref{eq:logistic}, showing the absolute error between the true and the recovered coefficient ${|r-\widehat{r}|}$ and delay ${|\tau-\widehat{\tau}|}$, together with RMSE$_{x^\prime}$ and RMSE$_{x}$ using brute force (BF), Bayesian optimization (BO) and particle swarm (PS).}\label{tab:logistic}
\end{sidewaystable}

\begin{figure}[!t]
\centering
\includegraphics[width=1\linewidth]{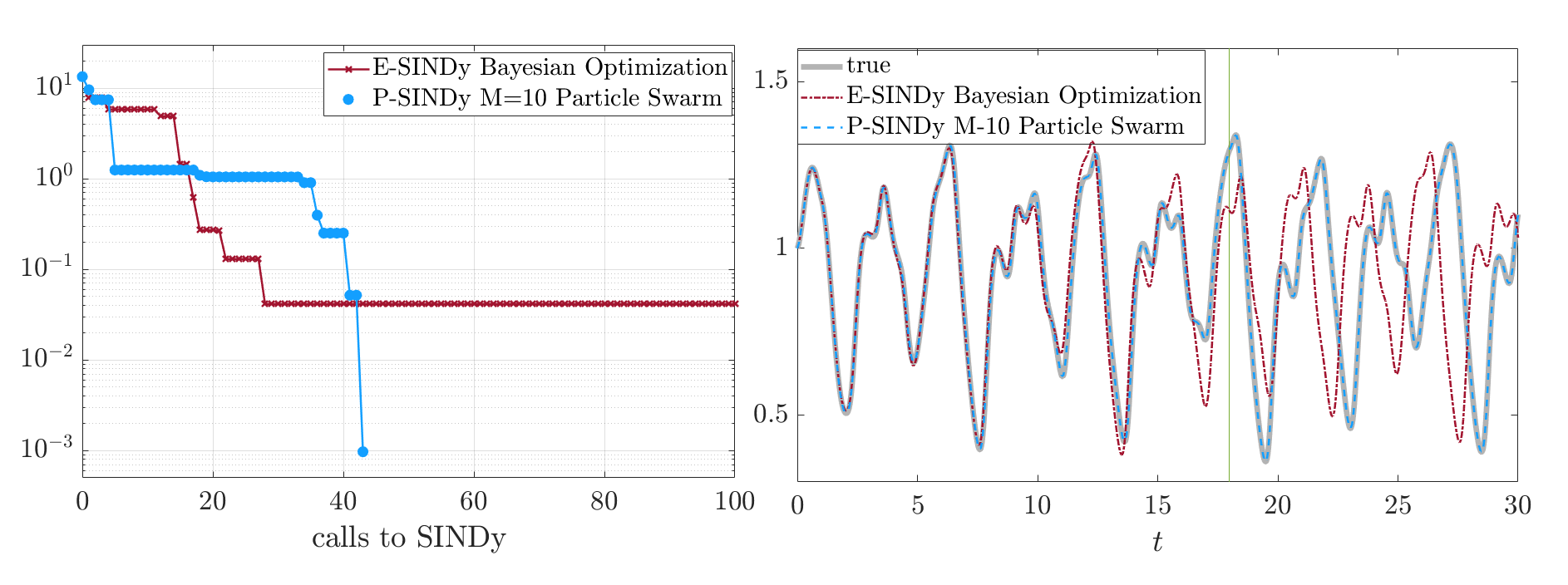}
\caption{(left) $\textrm{RMSE}_{x^\prime}$ for Mackey-Glass \eqref{eq:mg} for {E-SINDy} and {P-SINDy} with ${M=10}$ as a function of the total number of SINDy calls using Bayesian optimization (BO) and Paricle Swarm (PS). 
(right) Trajectory reconstruction for \eqref{eq:mg} showing the true trajectory (gray), {E-SINDy} trajectory using BO (red) and {P-SINDy} trajectory using ${M=10}$ and PS (blue). 
For the true model, the nominal values are ${(\tau,\alpha)=(1, 9.6)}$. 
{E-SINDy} with BO identify ${(\widehat{\tau},\widehat{\alpha})=(0.9996, 9.4969)}$, whereas {P-SINDy} with ${M=10}$ and PS identify ${(\widehat{\tau},\widehat{\alpha})=(1.0000, 9.6001)}$.}\label{fig:MG}
\end{figure}

A second series of tests concerns \eqref{eq:mg} with ${\beta=4}$ and ${\gamma=2}$, externally optimizing the unknown delay $\tau$ and Hill exponent $\alpha$ of nominal values 1 and 9.6, respectively. 
We use ${m=100}$ and ${T=30}$, with the initial function ${\varphi(s)=\cos(s)}$ for ${s\in [-1,0]}$, along with a polynomial library of degree ${d=2}$ that incorporates the Hill function from \eqref{eq:hill}. 
The representation of the RHS is depicted through the vector ${\Xi= [\,0,-\gamma,0,0,0,0,0,0,\beta,0\,]^{\top}}$ according to the library of candidate functions
\begin{equation}
    \Theta(X,X_{\tau},H)=
    \begin{bmatrix}
    1 & X & X_\tau & H & X^2 & XX_\tau & XH & X^2_\tau & X_\tau H & H^2
    \end{bmatrix},
\end{equation}
where $H$ is the matrix of samples of the Hill function \eqref{eq:hill}.
Figure~\ref{fig:MG} (left) compares the error $\textrm{RMSE}_{x^\prime}$ from BO and PS as a function of the total number of SINDy calls. 
In all three optimizers, the true values of $\tau$ and $\alpha$ are determined accurately, although the accuracy of BF again depends on the spacing of its candidate set, whereas BO and PS have other stopping criteria. 
Figure~\ref{fig:MG} (right) illustrates trajectory reconstructions using E-SINDy with BO and P10 with PS. 
The corresponding absolute error for these final optimized values along with RMSE errors appear in Table~\ref{tab:comparison2}.

\begin{figure}
\centering
\includegraphics[width=1\linewidth]{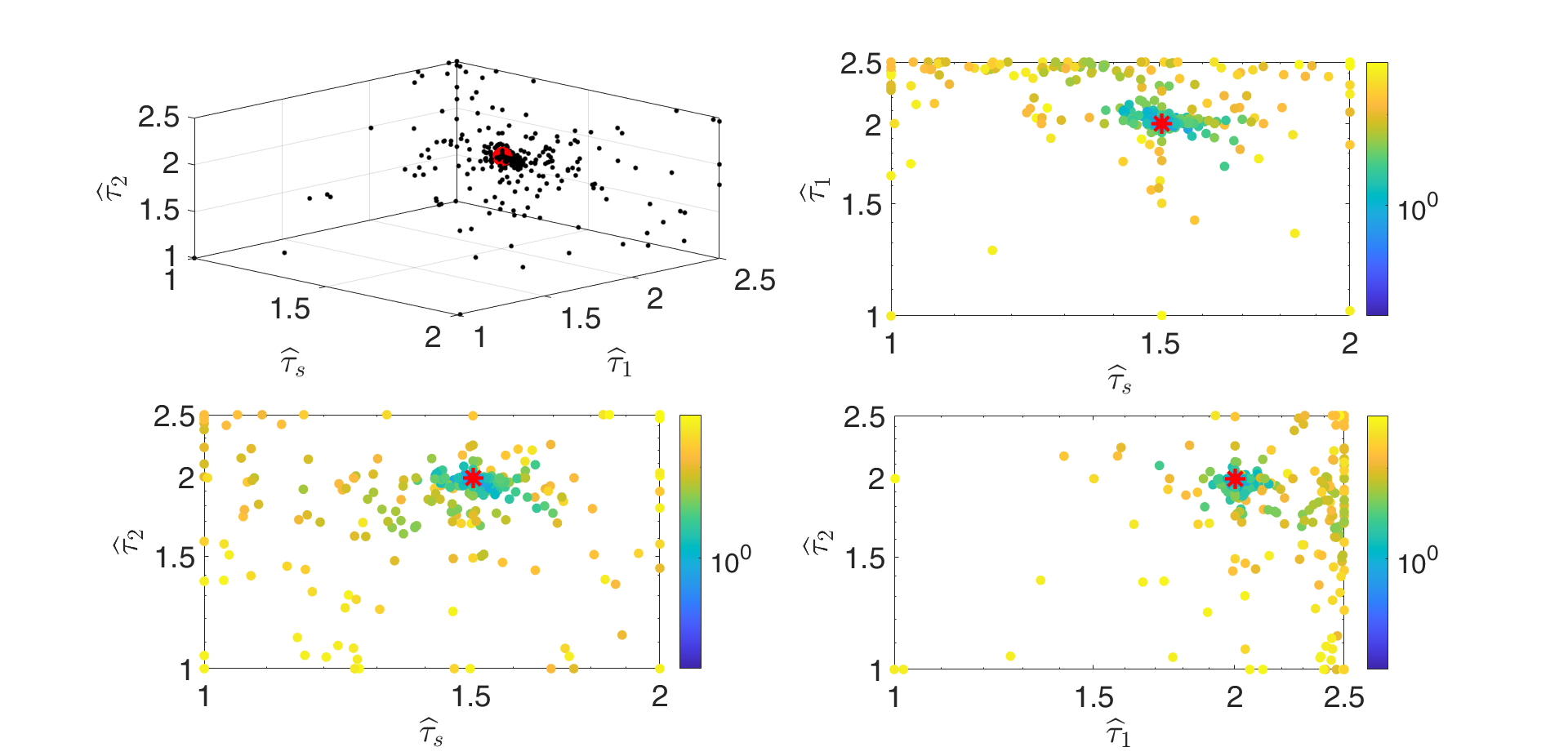}
\caption{Scatter plots of optimized unknown delays $(\widehat{\tau}_s,\,\widehat{\tau}_1,\,\widehat{\tau}_2)$ using particle swarm (PS) for the two-neuron model \eqref{eq:neuron}, color-coded by RMSE$_{x^\prime}$ values. 
Red stars mark the true values of the parameters.}\label{fig:neuron2} 
\end{figure}

Similarly, we collect data for the two-neuron model \eqref{eq:neuron}, setting the parameter values to ${\kappa = 0.5}$, ${\beta = -1}$, ${a_{12} = 1}$ and ${a_{21} = 2}$, with unknown delays of nominal values ${\tau_\mathrm{s} = 1.5}$, ${\tau_1 = 2}$ and ${\tau_2 = 2}$. 
Additionally, we choose ${m=100}$ and ${T= 30}$, along with a constant initial function $\varphi(s)\equiv [\,0.5,-0.5\,]^{\top}$, and utilize a polynomial library of degree ${d=2}$, incorporating trigonometric terms.  We optimize the multiple delays $\tau_\mathrm{s}$, $\tau_1$ and $\tau_2$ in the two-neuron model~\eqref{eq:neuron}. 
With {E-SINDy} all the delays are simultaneously optimized, whereas with {P-SINDy} only the maximum delay is optimized. 
Figure~\ref{fig:neuron2} provides an overview of these optimized delays using PS. 
Table~\ref{tab:comparison3} reports the absolute errors between recovered and true values and RMSE errors for BO and PS; we omit BF due to its large computational overhead and lower accuracy in multi-parameter searches.

Finally, Table~\ref{tab:comparison1} compiles the number of SINDy calls and CPU times for each approach (E, P5, P10, P15) combined with BF, BO or PS. 
In cases involving two or more unknown parameters (as in~\eqref{eq:mg} and~\eqref{eq:neuron}), the brute force optimizer become prohibitively expensive, leading us to omit BF from some experiments. 
When using {E-SINDy} with multiple delays, the overall cost scales with the product of each delay’s candidate set size. 
In contrast, {P-SINDy} only optimizes the maximum delay $\bar{\tau}$, considerably lowering the total computational load.

Increasing the collocation degree $M$ in {P-SINDy} does raise the cost; however, we often find that a low degree, such as ${M=5}$, already yields satisfactory results for many DDEs. 
The one exception here is the Mackey-Glass equation~\eqref{eq:mg}, whose potential chaotic nature demands ${M=10}$ to achieve accurate trajectory reconstruction. In that case, P10-PS converges in fewer SINDy evaluations, despite a similar CPU time.

\begin{table}[htbp]
\centering
\rotatebox{90}{%
\begin{minipage}{\textheight}
\centering

\small
\begin{tabular}{llrrrrrr}
\toprule
SINDy &optimizer& $|\xi_{2}-\widehat{\xi}_{2}|$ &$|\xi_{9}-\widehat{\xi}_{9}|$ & $|\tau-\widehat{\tau}|$ &$|\alpha-\widehat{\alpha}|$ &RMSE$_{x^\prime}$& RMSE$_{x}$\\
\midrule
 E &BO& \boldmath{$1.28\times^{-7}$} & \boldmath{$1.46\times10^{-7}$} & $4.00\times10^{-4}$ & $1.03\times10^{-1}$ & $7.74\times10^{-5}$ & \boldmath{$1.49\times10^{-3}$}\\
P5 &BO& $6.80\times10^{-5}$ & $4.12\times10^{-4}$ & $5.02\times10^{-2}$ & $3.67\times10^{-1}$ & $1.63\times10^{-4}$ & $9.85\times10^{-2}$\\
P10 &BO& $1.00\times10^{-6}$ & $1.30\times10^{-5}$ & \boldmath{$1.00\times10^{-4}$} & $1.34\times10^{-2}$ & $1.21\times10^{-5}$ & $1.64\times10^{-2}$ \\
P15 &BO& $7.19\times10^{-6}$ & $1.00\times10^{-6}$ & \boldmath{$1.00\times10^{-4}$} & \boldmath{$1.00\times10^{-4}$} & \boldmath{$5.45\times10^{-6}$} & $6.13\times10^{-3}$ \\ 
\midrule
 E &PS& $1.00\times10^{-6}$ & $2.00\times10^{-5}$ & $1.00\times10^{-6}$ & $1.00\times10^{-4}$ & $6.13\times10^{-7}$ & \boldmath{$5.45\times10^{-6}$}\\
P5 &PS& $2.06\times10^{-5}$ & $1.25\times10^{-6}$ & $6.00\times10^{-4}$ & $1.57\times10^{-1}$ & $1.64\times10^{-4}$ & $1.21\times10^{-2}$\\
P10 &PS& $1.03\times10^{-6}$ & $1.00\times10^{-6}$ & $2.30\times10^{-6}$ & $1.00\times10^{-4}$ & $1.09\times10^{-7}$ & $1.63\times10^{-4}$\\
P15 &PS& \boldmath{$1.00\times10^{-7}$} & \boldmath{$3.00\times10^{-7}$} & \boldmath{$1.80\times10^{-7}$} & \boldmath{$1.10\times10^{-5}$} & \boldmath{$1.03\times10^{-8}$} & $7.74\times10^{-5}$\\ 
\bottomrule
\end{tabular}
\caption{Performance comparison of {E-SINDy} (E) and {P-SINDy} (P, with ${M=5}$, ${10}$ and ${15}$) for \eqref{eq:mg}, highlighting the absolute error between the true and recovered coefficients ${|\xi_{i}-\widehat{\xi}_{i}|}$, delay ${|\tau-\widehat{\tau}|}$ and parameter ${|\alpha-\widehat{\alpha}|}$, together with $\textrm{RMSE}_{x^\prime}$ and $\textrm{RMSE}_{x}$ using Bayesian optimization (BO) and particle swarm (PS).}\label{tab:comparison2}

\vspace{0.5cm} 

\small
\begin{tabular}{llrrrrr}
\toprule
SINDy & optimizer &$|\tau_{s}-\widehat{\tau}_s|$ & $|\tau_{1}-\widehat{\tau}_1|$ & $|\tau_{2}-\widehat{\tau}_2|$ &RMSE$_{x^\prime}$& RMSE$_{x}$ \\ 
\midrule
 E &BO & $1.36\times10^{-2}$ & $5.22\times10^{-2}$ & \boldmath{$5.92\times10^{-2}$} & $7.74\times10^{-5}$ & $1.03\times10^{-1}$ \\
 P5 &BO  & - & - & $1.52\times10^{-1}$ & $1.63\times10^{-4}$ & $1.09\times10^{-1}$ \\
  P10 &BO & - & - & $1.06\times10^{-1}$ & $1.21\times10^{-5}$ & $2.47\times10^{-2}$ \\
 P15 &BO  & - & - & $8.27\times10^{-2}$ & \boldmath{$5.45\times10^{-6}$} & \boldmath{$1.88\times10^{-2}$}\\
 \midrule
 E  &PS & $2.70\times10^{-3}$ & $9.00\times10^{-4}$ & $4.90\times10^{-3}$ & $6.13\times10^{-4}$ & $2.76\times10^{-2}$ \\
 P5 &PS  & - & - & $3.60\times10^{-3}$ & $1.64\times10^{-6}$ & $6.21\times10^{-2}$ \\
  P10 &PS   & - & - & $1.90\times10^{-3}$ & $1.09\times10^{-7}$ & $4.85\times10^{-3}$ \\
 P15 &PS & - & - & \boldmath{$1.00\times10^{-6}$} & \boldmath{$1.03\times10^{-7}$} & \boldmath{$1.94\times10^{-3}$}\\
\bottomrule
\end{tabular}
\caption{Performance comparison of {E-SINDy} (E) and {P-SINDy} (P, with ${M=5}$, ${10}$ and ${15}$) for \eqref{eq:neuron}, highlightning the absolute error between the true and optimized delays ${|\tau_{s}-\widehat{\tau}_{s}|}$, ${|\tau_{1}-\widehat{\tau}_{1}|}$ and ${|\tau_{2}-\widehat{\tau}_{2}|}$, together with $\textrm{RMSE}_{x^\prime}$ and $\textrm{RMSE}_{x}$ using Bayesian optimization (BO) and particle swarm (PS).}\label{tab:comparison3}
\end{minipage}%
}
\end{table}

\begin{table}[htbp]
\centering
\small
\begin{tabular}{llrrrrrr}
\toprule
DDE & SINDy & \multicolumn{2}{c}{BF} & \multicolumn{2}{c}{BO} & \multicolumn{2}{c}{PS} \\
\cmidrule(r){3-4} \cmidrule(r){5-6} \cmidrule(r){7-8}
& & Calls & Time (s) & Calls & Time (s) & Calls & Time (s) \\
\midrule
\eqref{eq:logistic} & E & 1\,000 & 204 & 300 & 494& 308 & \textbf{67} \\
& P5 & 1\,000 & 217 & 300 & 438 & 241 & 70 \\
& P10 & 1\,000 & 449 & 300 & 485 & 179 & 100 \\
& P15 & 1\,000 & 1\,032 & 300 & 692 & \textbf{152} & 163 \\ \midrule
\eqref{eq:mg} & E & 10\,000 & hrs & 300 & 561 & 284 & \textbf{92} \\
& P5 & 10\,000 & hrs & 300 & 682 & 153 & 94 \\
& P10 & 10\,000 & hrs & 300 & 1\,080 & \textbf{51} & 102 \\
& P15 & 10\,000 & hrs & 300 & 1\,728 & 245 & 837 \\ \midrule
\eqref{eq:neuron} & E & - & - & 300 & 723 & 1\,303 & 906 \\
& P5 & - & - & 300 & 1\,081 & 344 & \textbf{378} \\
& P10 & - & - & 300 & 1\,410 & 342 & 494 \\
& P15 & - & - & 300 & 2\,088 & \textbf{130} & 580 \\
\bottomrule
\end{tabular}
\caption{Quantitative evaluation of computational performance of {E-SINDy} (E) and {P-SINDy} (P, with ${M=5}$, ${10}$ and ${15}$) showing number of SINDy calls and CPU times for the three different models \eqref{eq:logistic}, \eqref{eq:mg}, \eqref{eq:neuron} and for the optimizers brute force (BF), Bayesian optimization (BO) and particle swarm (PS).}\label{tab:comparison1}
 \end{table}

In conclusion, both {E-SINDy} and {P-SINDy} with the different external optimizers correctly recover the dynamics, identify unknown parameters and return accurate values for both $\textrm{RMSE}_{x'}$ and $\textrm{RMSE}_{x}$, with corresponding effective trajectory reconstruction for both training and testing. 
{E-SINDy} is more direct when the user already knows the structure and number of delays, whereas {P-SINDy} is more practical and in general computationally cheaper when multiple or unknown delays must be discovered from data. 
As far as the external optimizers are concerned, PS proves especially effective in minimizing the external search cost, often outperforming both BF and BO. As a general indication, P5 and P10 with PS represent on average the preferred choices.

\CCLsection{Neural Network Approaches}\label{sec:NN}

In this section, we introduce a set of approaches utilizing neural delay differential equations (NDDEs) to simultaneously learn the dynamics and the time delays from trajectory data. 
We first give a brief overview of modeling time delay systems with neural networks in Section~\ref{sec:NN_overview} and introduce the NDDE formalism in Section~\ref{sec:NN_formulate}.
Then we design the loss function and training algorithms in Section~\ref{sec:NN_train}.
We provide three examples of learning with NDDEs in Section~\ref{sec:NN_example} and discuss the size of data-driven models and compare the performance under different loss functions in Section~\ref{sec:NN_discuss}.

\CCLsubsection{Bridging Neural Networks with Dynamical Systems}
\label{sec:NN_overview}

\begin{figure}[!t]
\centering
\includegraphics[width=1\linewidth]{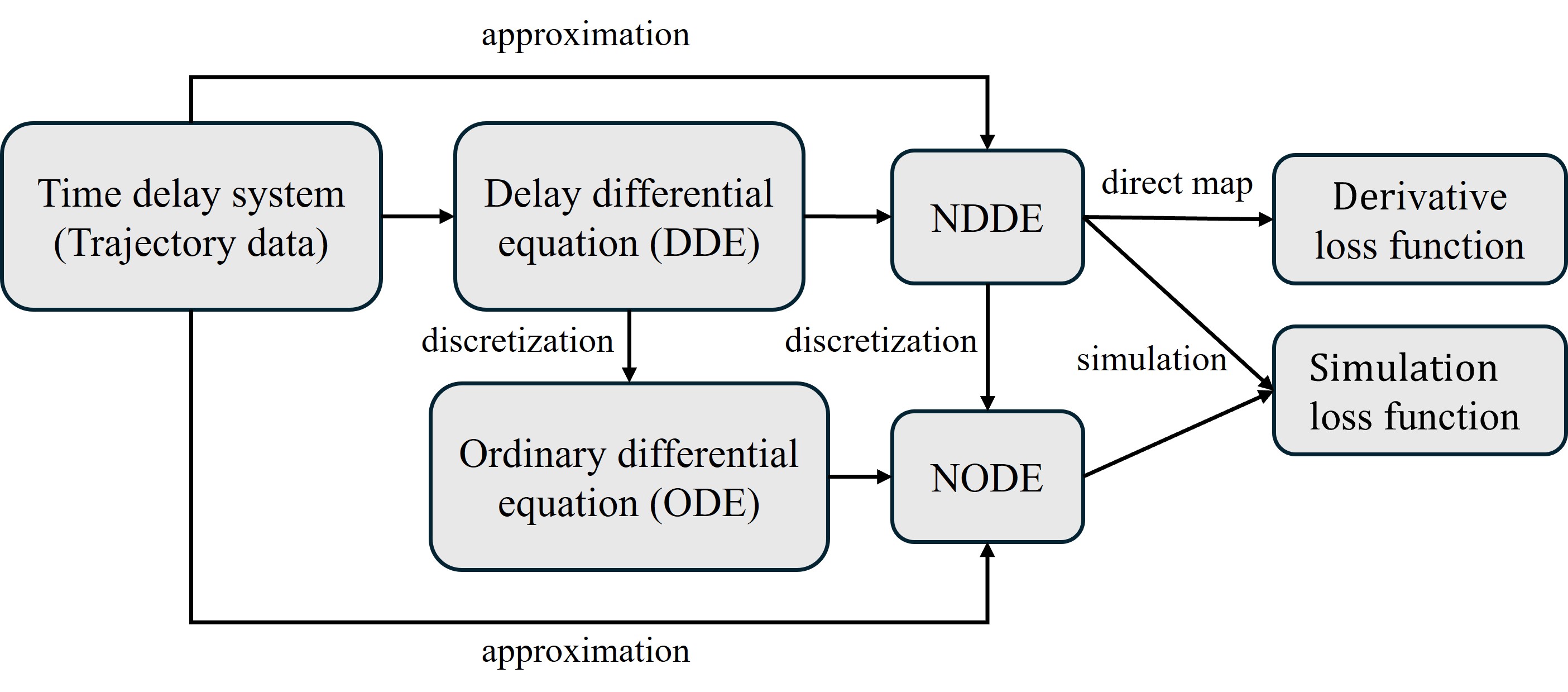}
\caption{Overview of modeling time delay systems with neural networks and differential equations}
\label{fig:NNmethod_overview}
\end{figure}

Among the data-driven methods developed for modeling dynamical systems, neural networks have shown their power in function approximation~\citep{pei2013mapping,bk19}. 
They are also in the midst of the most criticized methods, since they often involve very large number of parameters, which makes their analysis difficult~\citep{canziani2016analysis}. 
The generalizability of the trained network can be an issue if ``big data'' are not available, that is, neural networks can have poor performance outside the domain of the training data. 
Recent works made improvement by embedding prior knowledge---gained from first principle models---into the network design; see~\citep{ji2020,gupta2023neural,levine2022framework} for instance.
Adding time delays to neural networks as additional parameters can increase their capability of approximating the dynamics while keeping the network structure simple. 

In this section, we introduce a way to capture the dynamics of the time delay systems from data with the help of neural networks. 
We construct \emph{time delay neural networks with trainable delays} in continuous time to explicitly model the RHS of DDEs.

Figure~\ref{fig:NNmethod_overview} shows the overview of the learning framework.
As shown in Section~\ref{sec:DDEs}, time delay systems in continuous time can be described using DDEs, where the change rate of a dynamical system is determined by the past and current state of the systems; see \eqref{eq:ddef}.
Instead of deriving the relationship between the state derivative and the past states using first principles, one may use the flexible function approximators (i.e., time delay neural networks) to learn the relationship from the data directly.
As discussed for {P-SINDy} in Section~\ref{sec:PSINDy} it may be helpful to convert the DDEs into the ODEs via discretization.  
The same technique can be applied to convert the NDDEs into neural ordinary differential equations (NODEs), see the middle part of Figure~\ref{fig:NNmethod_overview}.

Below we introduce a framework of learning the time delay systems with neural networks and the algorithms can be categorized by the choice of loss functions, as shown on the right of Figure~\ref{fig:NNmethod_overview}.
If the loss function involves network simulations, one may simulate the NDDE directly or simulate a equivalent NODE (for which the training algorithms are well-developed and intensively studied).
If the loss function only contains a direct map between the data and the state derivatives, the conversion between NDDE and NODE is not necessary.
We will focus on the training algorithm using derivative loss and compare it to the simulation loss in one of the examples. 
References to other methods are provided for further explorations.

\CCLsubsection{Formulating the Neural Delay Differential Equations}
\label{sec:NN_formulate}

System~\eqref{eq:ddef} with $k$ discrete delays, ${0=\tau_0\leq \tau_1,\ldots,\tau_k \leq\tau_{\max}}$, can be re-written as
\begin{equation}\label{eq:DDE4}
    x^\prime(t) = f(\mathbf{x}_t),
\end{equation}
where ${\mathbf{x}_t = [\, x(t),x(t-\tau_1),\ldots,x(t-\tau_k)\, ]^{\top}\in{\mathbb{R}^{n(k+1)}}}$ and ${f:\mathbb{R}^{n(k+1)} \rightarrow \mathbb{R}^{n}}$. 
In order to capture the dynamics governed by~\eqref{eq:DDE4}, we construct the NDDE analogously:
\begin{equation}\label{eq:NDDE1}
    \hat{x}^\prime(t) = {\mathrm{net}}(\mathbf{\hat{x}}_t),
\end{equation}
with ${k+1}$ trainable delays, ${0\leq \hat{\tau}_0, \hat{\tau}_1,\ldots, \hat{\tau}_k \leq\tau_{\max}}$, and ${\mathrm{net}:\mathbb{R}^{n(k+1)}\rightarrow\mathbb{R}^n}$ representing the right hand side through a time delay neural network.
Suppose the network has $L$ layers with nonlinearity ${g_l(\cdot)}$ in each layer ${l=1,\ldots,L}$. 
Let us denote the input as ${\mathbf{z}_0^t}$ and the output as ${z_L^t}$ at time $t$. 
Then we have
\begin{equation}\label{eq:network}
z^t_{L} ={\mathrm{net}}(\mathbf{z}^t_{0}),
\end{equation}
where
\begin{equation}\label{eq:NDDE2}
    \begin{split}
    z^t_{1}&=g_1(W_{1}\mathbf{z}^t_{0} +b_{1}),
    \\
    z^t_l &= g_l\big(W_l z^t_{l-1} + b_l\big), \quad l=2,\ldots,L,
    \end{split}
\end{equation}
where $W_l$ and $b_l$ denote the weights and biases, respectively, and one may use different nonlinearities $g_l$ including the rectified linear unit (ReLU) and different sigmoidal functions.
This network can be used to approximate the state derivative of the delayed dynamical system by defining ${z^t_L=\hat{x}^\prime(t)}$ and ${\mathbf{z}^t_{0} =\mathbf{\hat{x}}_t}$, cf.~\eqref{eq:NDDE1} and \eqref{eq:network}.

\begin{figure}[!t]
      \centering
    \includegraphics[width= 10cm]{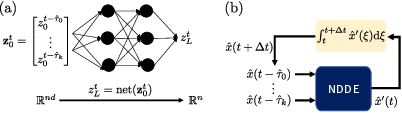}
   \caption{(a) Illustration of a time delay neural network constructed to capture the dynamics of a neural delay differential equation (NDDE) with a single delay. 
   (b) Illustration of the simulation of an NDDE.}
      \label{fig:NN}
\end{figure}

Note that some of the delays trained in an NDDE do not have to be in an ascending/descending order and they can potentially have the same value. 
When the NDDE \eqref{eq:NDDE1} is used to approximate the dynamics of the time delay system~\eqref{eq:DDE4}, 
knowing the number of delays in the system is beneficial but not required~\citep{ji2024nd}.
Moreover, if the existence of a non-delayed term with ${\tau_0=0}$ is known, one may fix the delay value ${\hat{\tau}_0 = 0}$ to further restrict the neural network.
Figure~\ref{fig:NN}(a) illustrates a time delay neural network used to capture the right hand side of an NDDE with discrete delays, while panel (b) shows a block diagram illustrating the corresponding NDDE simulation.
We emphasize that below we construct time delay neural networks with trainable delays, that is, apart from the weights $W_l$ and biases $b_l$, ${l=1,\ldots,L}$ in \eqref{eq:NDDE2}, we also aim to learn the delays ${\hat{\tau}_0, \hat{\tau}_1,\ldots, \hat{\tau}_k}$ themselves.

For example, an autonomous delay dynamical system with one delay (${k = 1}$) can be described by
\begin{equation}
\label{eq:DDE3}
    \dot{x}(t) = f(x(t),x(t-\tau)),
\end{equation}
whose RHS $f$ and delay $\tau$ may be unknown.
Our goal is to learn $f$ as well as $\tau$ from the data using the NDDE
\begin{equation}
\label{eq:NDDE3}
    \hat{x}^\prime(t) = {\mathrm{net}}(\hat{x}(t),\hat{x}(t-\hat{\tau})),
\end{equation}
with trainable delay ${\hat{\tau}}$.

Note that control systems with state feedback can also be viewed as autonomous systems.
For example, with slight abuse of notation, the system 
${\dot{x}(t) = f(x(t),u(t))}$ with delayed state feedback ${u(t) = k(x(t-\tau))}$ yields
\begin{equation}
    x^\prime(t) = f\big(x(t), k(x(t-\tau))\big).
\end{equation}
When the nonlinear dynamics $f$ of the system, the feedback law $k$, and the delay $\tau$ are not known, we need to learn these from trajectory data.

The trained neural network~\eqref{eq:NDDE3} can then be solved by any DDE solver for a given initial history.
The state at time $t$ is given by
\begin{equation}
\label{eq:DDEsolver}
    \hat{x}(t) = {\mathrm{DDEsolver}}({\mathrm{net}},x_{t_0},t),
\end{equation}
where the history function ${x_{t_0}(s), s\in [t_0-\tau_{\max}, t_0]}$ is approximated from data ${x(t_0),x(t_0-\Delta t),\ldots,x(t_0-\tau_{\max})}$ and ${\Delta t}$ denotes the sampling time.

Due to the trainable delays in the networks and the recurrent feature of network simulation, the time delay neural network is able to take the time dependency into consideration, while keeping the number of parameters relatively low.
Moreover, once the dynamics of the DDE~\eqref{eq:NDDE1} are captured by the neural network~\eqref{eq:network}, they can be analyzed using tools from the dynamical systems and control literature.
Additionally, with explicit delay parameters in the input layer, the network is more interpretable and generalizable.

\CCLsubsection{Design the loss function and training algorithms}
\label{sec:NN_train}

In this section, we construct loss functions and develop training algorithms for the neural networks illustrated in Figure~\ref{fig:NN}(a). 
We categorize the training algorithms based on whether the simulation of the NDDE is involved in the training process, which is also related to the type of loss function used, see Figure~\ref{fig:NNmethod_overview}.

First, let us consider the derivative loss, where the network is not simulated during training and it predicts the state derivative directly from the data.
This corresponds to substituting the NDDE \eqref{eq:NDDE1} with 
\begin{equation}\label{eq:NDDE1_mod}
    \hat{x}^\prime(t) = {\mathrm{net}}(\mathbf{x}_t),
\end{equation}
where the data enters the right hand side directly through $\mathbf{x}_t$ (notice the lack of hat).
For  example, in case of the autonomous system~\eqref{eq:NDDE3} with one delay, such direct map is given by
\begin{equation}
\label{eq:direct}
    \hat{x}^\prime(t) = {\mathrm{net}}(x(t),x(t-\hat{\tau})),
\end{equation}
where ${x(t),x(t-\hat{\tau})}$ are obtained from the data and ${\hat{\tau}}$ is the learned delay at the current iteration. 
Since ${\hat{\tau}}$ is a continuous variable, ${x(t-\hat{\tau})}$, i.e., the value of the state $x$ at time ${t-\hat{\tau}}$, may not exist in the discrete measurement data. 
Therefore, linear interpolation is used to recover the delayed state ${x(t-\hat{\tau})}$ from the discrete history data.

The derivative loss will take the form
\begin{equation}
\label{eq:loss_deriv}
    \mathcal{L}_{x^\prime} =\frac{1}{n N} \sum_{m \in \mathrm{batch}} \sum_{j = 1}^{n}w_m\Big(\hat{x}^\prime_j(t_m)-x^\prime_j(t_m)\Big)^2,
\end{equation}
where ${t_m = m\Delta t}$, is the sampling time, $N$ is the number of samples used in one update, $w_m$ represents the loss weight associated with the $m$-th data point, and ${x_{j}}$ denotes the $j$-th component of the vector ${x \in \mathbb{R}^n}$.
It is shown that emphasizing the initial segment of the data helps to learn the dynamics. In the examples of this chapter, we are using ${w_m = 1/t_m}$, same as in~\citep{ji2023l4dc}.
Since the derivative loss does not require solving the DDE, the ${\hat{x}^\prime}$ is given by~\eqref{eq:NDDE1_mod}.

As a comparison, we also introduce the simulation loss function to measure the error between the NDDE simulation and the data:
\begin{equation}\label{eq:loss_simu}
    \mathcal{L}_{x} = \frac{1}{nHN}\sum_{m \in \mathrm{batch}} \sum_{i = 1}^{H} \sum_{j = 1}^{n} \Big(\hat{x}_j(t_m+i\Delta t)-x_j(t_m+i\Delta t)\Big)^2,
\end{equation}
where ${H\Delta t}$ is the simulation horizon. 
As the data is sampled at ${t = t_m+i\Delta t}$, the network simulation is also evaluated at the same time moments. 
Here, $\hat{x}$ is given by~\eqref{eq:DDEsolver}.

Both loss functions~\eqref{eq:loss_deriv} and~\eqref{eq:loss_simu} depend on the parameters in the time delay neural networks, including the weights $W_l$ and biases $b_l$, ${l = 1,\ldots,L}$, and the delays $\hat\tau_r$, ${r = 0, 1,\ldots,k}$.
Any gradient-based algorithms require the calculation of the gradient with respect to each parameter in the set ${\{W_l,b_l,\hat{\tau}_r\}}$, that is, ${\frac{\partial \mathcal{L}}{\partial W_l}}$, ${\frac{\partial \mathcal{L}}{\partial b_l}}$, and ${\frac{\partial \mathcal{L}}{\partial \hat{\tau}_r}}$. 

The gradient of the loss with respect to delay $\hat{\tau}_r$ is given by
\begin{equation*}
    \frac{\partial \mathcal{L}}{\partial \hat{\tau}_r} = \frac{\partial \mathcal{L}}{\partial  x(-\hat{\tau}_r)} \frac{\partial x(-\hat{\tau}_r)}{\partial \hat{\tau}_r} = -\frac{\partial \mathcal{L}}{\partial x(-\hat{\tau}_r)} x^\prime(-\hat{\tau}_r).
\end{equation*}
Here $\hat{\tau}_r$ is a continuous variable, so ${x(-\hat{\tau}_r)}$ and ${x^\prime(-\hat{\tau}_r)}$ are approximated from the linear interpolation of the history samples.
Namely, one can express any delay value as ${\tau = (i+\alpha)\Delta t}$ with ${i\in\mathbb{N}}$ and ${\alpha \in [0,1)}$, and approximate the state at ${t-\tau}$ as
\begin{equation*}
\label{eq:xtau}
    x(t-\tau)= x(t-(i+\alpha)\Delta t) \approx(1-\alpha)x(t-i\Delta t)+\alpha x(t-(i+1)\Delta t).
\end{equation*}
The state derivatives can be obtained as follow
\begin{equation*}
\begin{split}
    \dot{x}(t-\tau) 
    &\approx \frac{x(t-i\Delta t)-x(t-(i+\alpha)\Delta t)}{\alpha\Delta t} \\
    &\stackrel{(\ref{eq:xtau})}{\approx}\frac{x(t-i\Delta t)-x(t-(i+1)\Delta t)}{\Delta t},
\end{split}
\end{equation*}
where Euler's method is used in the first approximation.

\begin{algorithm} [tp]\label{alg:training}
\begin{minipage}{0.88\linewidth}
\textbf{Goal:} Learn ${p = \{W_l,b_l, \hat{\tau}_r \}, l = 1,\ldots,L, r = 0,1,\ldots, k}$ from data.
{\begin{itemize}
    \item Choose the hyper parameters according to the loss function. Choose the learning rates $\eta_\tau$ and $\eta$ for delays and other parameters.
    \item Set the maximum iteration number $q_{\mathrm{max}}$
    and the maximum allowed delay $\tau_{\mathrm{max}}$.
    \item Initialize $p_0$, e.g., initialize $\hat{\tau}_r$ uniformly in the interval ${[0,\tau_{\mathrm{max}}]}$, $b_l$ as zeros and $W_l$ using Glorot initialization~\citep{glorot2010understanding}.
 \end{itemize}
}
 \For{${q = 1,\ldots, q_{\mathrm{max}}}$}{
 \begin{minipage}{0.9\linewidth}
  \begin{itemize}
    \item Randomly take $N$ segments or samples from the training 
    data. 
    \item Calculate the loss $\mathcal{L}_q$.
    \item Back-propagate to get the gradients with respect 
    to parameters using \texttt{dlgradient}.
    \item Update $p_q$ from $p_{q-1}$ using \texttt{adamupdate} and
    impose positivity constraint on the delays $\hat{\tau}_r$.
 \end{itemize}
   
   \If{$ q > 1$ {\textbf{and}} $\mathcal{L}_{q} < \mathcal{L}_{q-1}$}{
  $p_{\mathrm{best}} =p_{q}$}
 \end{minipage}
 }{$p = p_{\mathrm{best}}$, evaluate the NDDE on testing datasets}
 \caption{Training algorithm for the time delay neural network.}
\vspace{3mm}
 \label{alg:NDDE_learning}
\end{minipage}
\end{algorithm}

After the gradients ${\{\frac{\partial \mathcal{L}}{\partial W_l}, \frac{\partial \mathcal{L}}{\partial b_l}, \frac{\partial \mathcal{L}}{\partial \hat{\tau}_r}\}}$ are calculated via backpropagation, one can update the parameters using any gradient-based methods.
Here we provide the adaptive moment estimation (ADAM) algorithm devised by~\cite{kingma2014adam}), and the examples in the following sections use the \texttt{adamupdate} function from MATLAB Deep Learning Toolbox.

The ADAM formula for updating each individual parameter is given by
\begin{equation*}
    p_{q+1} = p_q-\frac{\eta}{\sqrt{\hat{v}_q}+\epsilon}\hat{m}_q
\end{equation*}
with 
\begin{equation*}
\hat{m}_q = \frac{m_q}{1-\beta_1^q},\qquad \hat{v}_q= \frac{v_q}{1-\beta_2^q},
\end{equation*}
where $m_q$ and $v_q$ are the first moment estimate and second moment estimate at iteration $q$, respectively.
The updates for $m_q$ and $v_q$ are given by
\begin{equation*}
    m_q = \beta_1 m_{q-1} + (1-\beta_1) g_{q},\qquad    v_q =\beta_2 v_{q-1} + (1-\beta_2) g_{q}^2,
\end{equation*}
using the gradient information $g_{q}$ at iteration $q$. The parameters ${\eta,\epsilon,\beta_1,\beta_2}$ can be tuned by the user to achieve better performance; see~\citep{kingma2014adam}.

In addition to the parameter updates with the gradient information, we restrict the delay value to be between $0$ and $\tau_{\mathrm{max}}$ considering the positivity of the delay and the length of memory. 
An example of the NDDE training algorithm is illustrated in Algorithm~\ref{alg:NDDE_learning}.
One may customize the training algorithm by designing different loss functions, selecting parameter update rules and the hyper parameters, and implement parallel training techniques. 
We provide the MATLAB codes for all the examples on \href{https://github.com/Jxb814/ch9_NDDE_sample_code/tree/main}{github}.

\CCLsubsection{Examples and Results}
\label{sec:NN_example}

In this section, we test the proposed learning framework with time delay neural networks on multiple systems.
The examples presented in this section are all scalar systems, with one or two delays.
An additional example of learning multiple-state autonomous system will be given and compared to SINDy approach in Section~\ref{sec:compare}.
For the delay logistic equation and the Mackey-Glass equation, the NDDE with one trainable delay takes the form
\begin{equation}
\label{eq:NDDE_auto}
    \hat{x}^\prime(t) =
  W_3 \tanh{\Big(W_2 \tanh{\big(W_1 
  \begin{bmatrix}
  \hat{x}(t)
  \\  
  \hat{x}(t-\hat{\tau})
  \end{bmatrix}
  +b_1\big)+b_2\Big)}},
\end{equation}
where ${W_1 \in \mathbb{R}^{l_{1}\times 2}}$, ${W_2 \in \mathbb{R}^{l_{2}\times l_{1}}}$, ${W_3 \in \mathbb{R}^{1\times l_{2}}}$, ${b_1 \in \mathbb{R}^{l_{1}\times 1}}$, ${b_2 \in \mathbb{R}^{l_{2}\times 1}}$ and $l_1,l_2$ are the number of neurons in each hidden layer.
We emphasize that in the training process with the derivative loss (when there is no simulation), the actual input of the NDDE is ${[\, x(t), x(t-\hat{\tau})\, ]^\top}$ from the data, instead of  ${[\, \hat{x}(t), \hat{x}(t-\hat{\tau})\, ]^\top}$ from the network simulation, cf.~\eqref{eq:direct} and \eqref{eq:NDDE1_mod}.

\CCLsubsubsection{Delay logistic equation}

\begin{figure}[!t]
    \centering
    \includegraphics[width=1\linewidth]{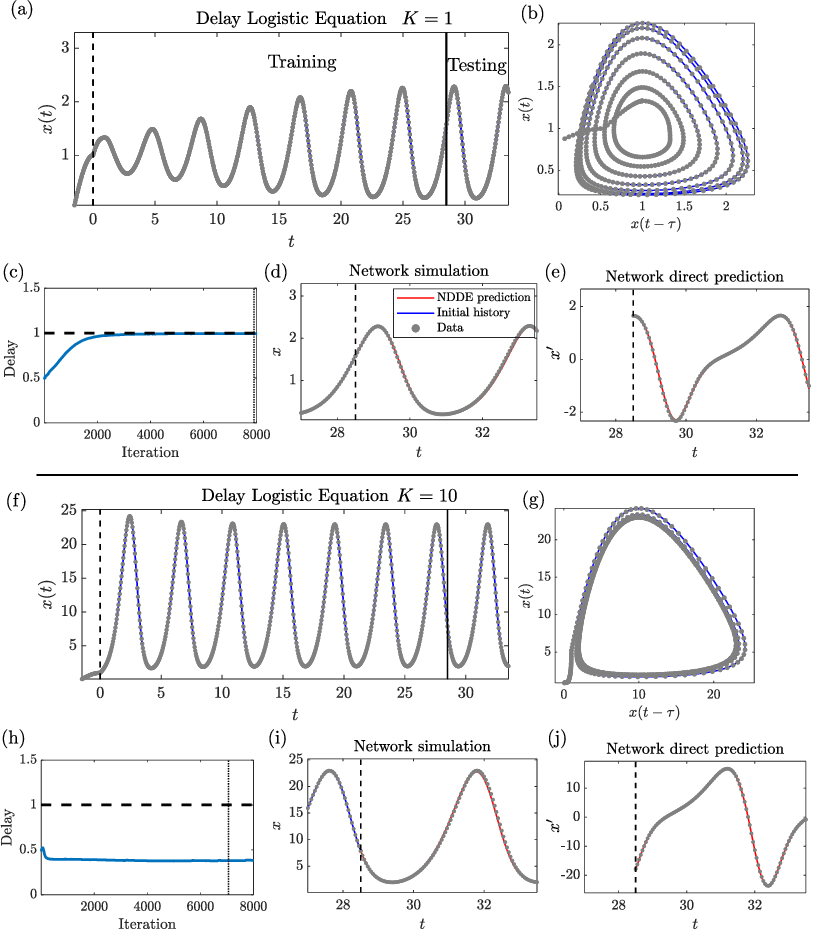}
    \caption{Training NDDE \eqref{eq:NDDE_auto} (2 hidden layers and 10 neurons in each layer) with two different datasets of delay logistic equation.}
    \label{fig:delay_logistic_data}
\end{figure}

First, let us consider the delay logistic equation \eqref{eq:logistic}.
By choosing ${K = 1}$, ${\tau = 1}$ and initial history ${\varphi(s) = \cos (s)}$, we generate simulation data with a sampling time ${\Delta t = 0.5}$. 
The first section of length ${T_{\mathrm{tr}} = 30}$ (including the initial history) is used as training data and the next section of length ${T_{\mathrm{st}}=5}$ is used as testing data to examine the performance of trained model. 
The data is shown in Figure~\ref{fig:delay_logistic_data}(a) as a function of time and the training data is also demonstrated in panel (b) in (projected) state space.
The initial delay of the NDDE is set to be 0.5 and it converges to the true value ${\tau = 1}$ during training, see panel (c).
The maximum training iteration is ${q_{\max} = 8000}$ for this example and the dotted vertical line indicates the iteration where the training loss is minimal (i.e., ${p=p_{\mathrm{best}}}$).
Panel (d) shows the simulations of the trained NDDE \eqref{eq:NDDE_auto} on the testing dataset, i.e., only the data before the dashed vertical line is provided as history in the simulation.
The direct prediction (with ${[\, x(t), x(t-\hat{\tau})\, ]^\top}$ in the RHS of the trained NDDE \eqref{eq:NDDE_auto}) is shown in Figure~\ref{fig:delay_logistic_data}(e) for the testing data.
Note that the direct prediction does not solve the differential equation, therefore the NDDE prediction in panel (e) is not the numerical differentiation of the simulation result in panel (d).

When the dynamics of the data is not rich enough, the NDDE can possibly learn a ``wrong'' delay while still recover the data in both simulation and direct derivative prediction. 
To demonstrate this, we generate another trajectory (of the same length) using ${K = 10}$, ${\tau = 1}$, where the training data is concentrated around a limit cycle.
The training data is presented in Figure~\ref{fig:delay_logistic_data}(f) and (g) and the convergence of the learned delay is shown in panel (h).
Without any prior knowledge or restrictions on the network design, the delay ${\tau = 1}$ is not learned correctly. 
Nevertheless, the network simulation and direct prediction are both close to the testing data, as shown in panels (i) and (j). This means that both simulation and derivative loss surface has local minima at this delay. From the point of view of approximating the data, the algorithm delivers good results while the truth about the delay still remains hidden.

\CCLsubsubsection{Mackey-Glass equation}

\begin{figure}
    \centering
    \includegraphics[width=1\linewidth]{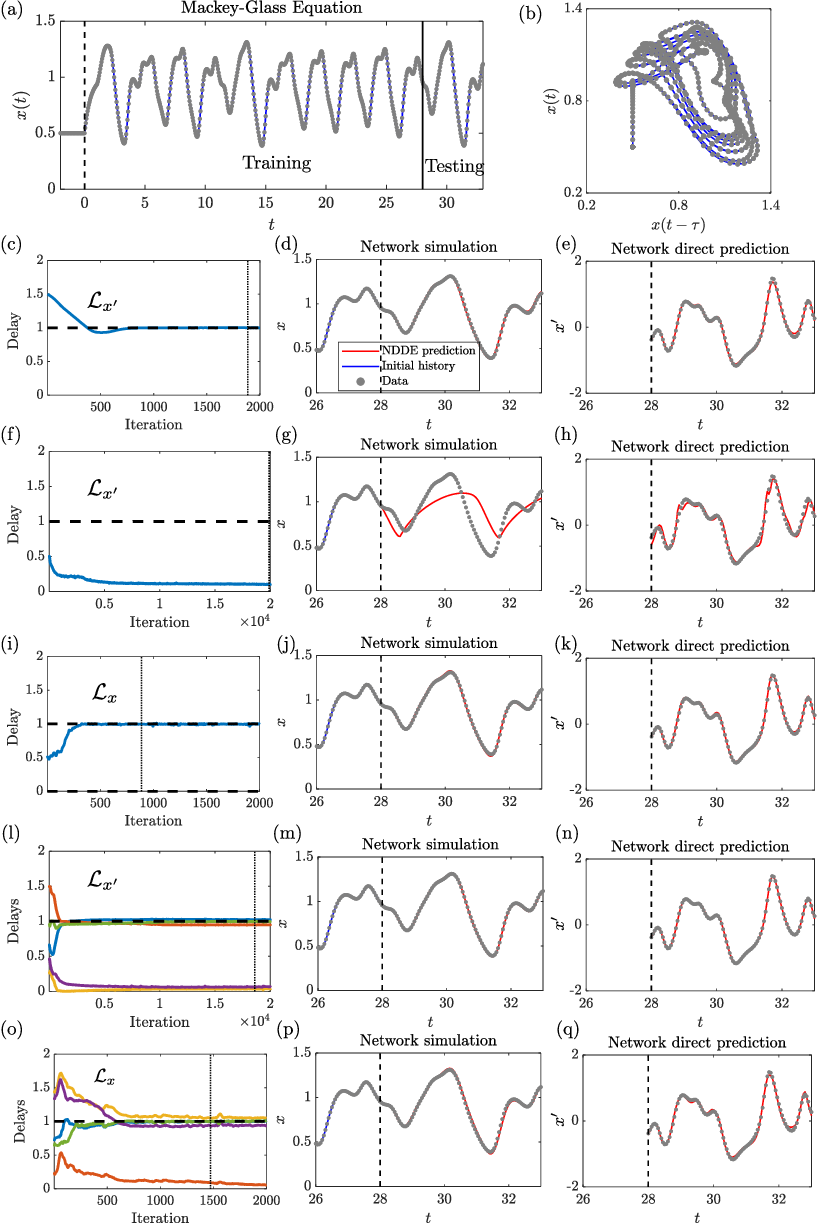}
    \caption{Training NDDEs \eqref{eq:NDDE_auto} and \eqref{eq:NDDE_auto5} (2 hidden layers, 5 neurons in each) on simulation data generated by Mackey-Glass equation.}
    \label{fig:mg_data}
\end{figure}

Now we test the same algorithm on the Mackey-Glass equation \eqref{eq:mg} with parameters ${\alpha = 9.6}$, ${\beta = 4}$, ${ \gamma = 2}$ and delay ${\tau = 1}$.
We generate simulation data with a sampling time ${\Delta t = 0.5}$ from a constant history ${\varphi(s)\equiv0.5}$ shown in Figure~\ref{fig:mg_data}(a), and the phase portrait of the training data are depicted in panel (b).
This system is chaotic under the given parameters resulting in a rich training data. 

Here we demonstrate the performance of the trained NDDE \eqref{eq:NDDE_auto} on the testing data, considering different initial delay values, loss functions and network structures. 
If the derivative loss function is used in the training process, the delay may not converge to our desired value when the initial value of the delay is small. 
Figure~\ref{fig:mg_data}(c) and (f) show the evolution of the trainable delay for  initial delays 1.5 and 0.5, respectively. 
When the initial delay is 1.5, the delay can be learned through 2000 iterations and the trained NDDE is able to give good predictions for the testing data, both in simulation and in direct derivative predictions, see panels (d) and (e).
When the initial delay is 0.5, the trainable delay converges to a small, incorrect  value. 
In this case, while the direct derivative prediction of the trained network performs relatively well for the testing data as shown in panel (h), the network simulation deviates significantly from the data in panel (g). 
This indicates an undesired local minimum of the derivative loss function around a small delay value.

We remark that it is possible to learn numerical differentiation schemes when training NDDEs because those schemes indeed constitute local minima of the derivative loss function, see the study in~\citep{ji2021l4dc}.
For example, the NDDE can approximate
\begin{equation}
    \hat{x}^\prime = \frac{1}{\hat{\tau}}x(t)-\frac{1}{\hat{\tau}}x(t-\hat{\tau}),
\end{equation}
instead of the right hand side of the differential equation.
This can happen when data are sampled densely (the undesired solution may be the global minimum) and when the initial delay is small (starts too close to the undesired solution).
To address this issue, one may implement parallel training (i.e., train the network from multiple initial delays) or choose a loss function which includes multi-step simulation.
For instance, the NDDE with initial delay 0.5 converges to the true value under the simulation loss \eqref{eq:loss_simu} with horizon ${H = 10}$ as shown in Figure~\ref{fig:mg_data}(i) and the trained NDDE performs well on the testing data.
While using the simulation loss function one can avoid learning a small delay, the computation time increases significantly due to the simulation involved in the training process. 
Comparison between the computation time of two loss functions will be presented in the next section.

Besides the NDDE~\eqref{eq:NDDE_auto}, we also test the training algorithm with multiple trainable delays, assuming no knowledge about the number of delays in the system.
The last two rows of Figure~\ref{fig:mg_data} shows the results for the NDDE
\begin{equation}
\label{eq:NDDE_auto5}
    \hat{x}^\prime(t) =
  W_3 \tanh{\Big(W_2 \tanh{\big(W_1 
  \begin{bmatrix}
  \hat{x}(t)
  \\  
  \hat{x}(t-\hat{\tau}_1)
  \\  
  \hat{x}(t-\hat{\tau}_2)
  \\  
  \hat{x}(t-\hat{\tau}_3)
  \\  
  \hat{x}(t-\hat{\tau}_4)
  \\  
  \hat{x}(t-\hat{\tau}_5)
  \end{bmatrix}
  +b_1\big)+b_2\Big)}},
\end{equation}
with 5 trainable delays where the dimensions of the weight matrix ${W_1 \in \mathbb{R}^{l_{1}\times 5}}$ is adjusted.
Both the derivative and the simulation loss can be used to learn the delay---the trained delays are concentrated around the true value ${\tau=1}$ and zero, the latter corresponding to the non-delayed term in \eqref{eq:mg}. 
In both cases the trained NDDE predicts the testing data well.

\CCLsubsubsection{Climate model}

In addition to nonlinear autonomous systems with single delay, we can extend the training algorithms to systems with multiple delays and external forcing. 
Here we use an NDDE with similar structure to that of \eqref{eq:NDDE_auto} in order to learn a climate model from trajectory data.

\begin{figure}[!t]
    \centering
    \includegraphics[width=0.8\linewidth]{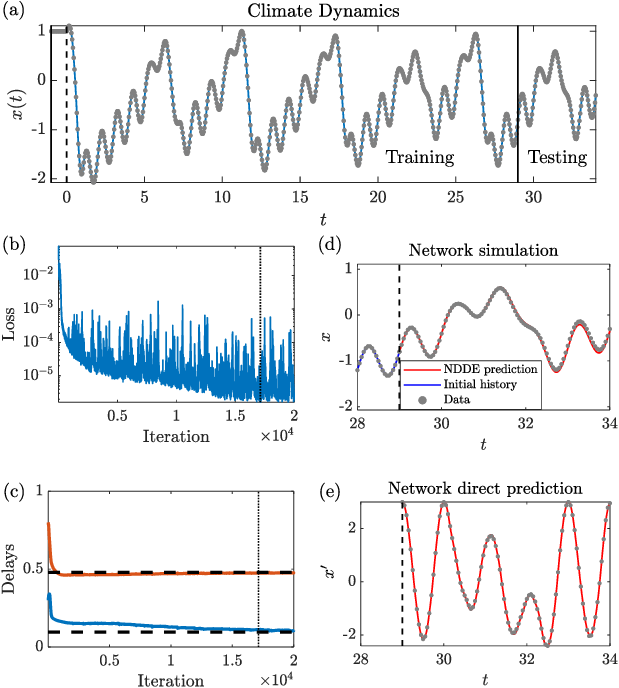}
    \caption{Training and testing performance of the climate model. The NDDE~\eqref{eq:NDDE_climate} has 10 neurons in each hidden layer. 
    }
    \label{fig:climate_data}
\end{figure}

\begin{example}\label{ex:climate}
We consider a \textit{climate model} with dynamics 
\begin{align}\label{eq:climate_model}
    x^\prime(t) = aA(\kappa,x(t-\tau_1))-bA(\kappa,x(t-\tau_2))+c u(t),
    \\
    A(\kappa,x) = 
    \begin{cases}
        d_u \tanh (\frac{\kappa}{d_u}x), & \text{if $x\geq0$},
        \\
        d_l \tanh (\frac{\kappa}{d_l}x), & \text{if $x<0$}.
    \end{cases}
\end{align}
This model describes the deviation of thermocline depth from its long-term mean at the eastern boundary of the equatorial Pacific Ocean. 
The chaotic behavior of this DDE is used to demonstrate the irregular phenomenon of the El Ni\~{n}o Southern Oscillation (ENSO). 
More detailed discussion can be found in~\citep{saunders2001boolean,keane2017climate,keane2019effect}.
\end{example}

In this example, the positive and negative feedbacks are associated with two different fixed delays 
${\tau_1 = 0.0958}$ and ${\tau_2 = 0.4792}$, 
while the other parameters are 
${a = 2.02}$, ${b = 3.03}$, ${c = 2.6377}$, ${d_u = 2.0}$, ${d_l = -0.4}$. 
We generate a trajectory of length ${T = 30}$  using the constant history ${\varphi(s)\equiv 1}$ and a periodic input ${u(t) = \cos (2\pi t)}$ which represents the seasonal forcing.
The data is shown in Figure~\ref{fig:climate_data}(a) with a sampling time ${\Delta t= 0.05}$.

For this learning task, we use a time delay neural network which takes the forcing $u(t)$ and the state history $x_t$ as input and gives the state derivative $x^\prime(t)$ as output. 
The network consists of 2 hidden layers with 5 neurons in each layer.
With two trainable delays, the NDDE takes the form
\begin{equation}
\label{eq:NDDE_climate}
    \hat{x}^\prime(t) =
  W_3 \tanh{\Big(W_2 \tanh{\big(W_1 
  \begin{bmatrix}
  \hat{x}(t-\hat{\tau}_1)
  \\  
  \hat{x}(t-\hat{\tau}_2)
  \\
  u(t)
  \end{bmatrix}
  +b_1\big)+b_2\Big)}},
\end{equation}
where $u(t)$ is obtained from interpolating the recorded external forcing data.

The training loss and the delays along the iterations are shown in Figure~\ref{fig:climate_data}(b) and (c), while the network simulation and the direct derivative prediction of the trained NDDE are presented in panels (d) and (e). 
The NDDE is able to learn both delays and predict the short-term future.

\CCLsubsection{Evaluation and Discussion}
\label{sec:NN_discuss}

The performance of the trained NDDEs on the testing datasets are shown in Tables~\ref{tab:logistic_NDDE},~\ref{tab:mg_NDDE} and~\ref{tab:climate_NDDE} for the three examples above in terms of the learned delays and the RMSE. 
We applied the same training algorithm \ref{alg:training} with derivative loss \eqref{eq:loss_deriv} using fully-connected 2-layer neural networks of different sizes. 
As the number of neurons in the hidden layers are increased, the performance of the trained networks improves. 
Large improvements are observed when increasing the number of neurons from 2 to 5, beyond which the improvements are not significant.

\begin{table}[!t]
    \centering
    \small
    \begin{tabular}{ccccc}
    \toprule
    Loss & Size  & RMSE$_{x^\prime}$ & RMSE$_{x}$ &$|\tau-\hat{\tau}|$ \\
       \midrule
       \multirow{4}{*}{$\mathcal{L}_{x^\prime}$} & 2 & $3.42\times10^{-1}$  & $1.83\times10^{-1}$ & $1.70 \times 10^{-1}$\\
      & 5 & $6.19\times10^{-2}$ & $3.38\times10^{-2}$ & $1.76 \times 10^{-2}$\\
       & 10 & $2.43\times10^{-2}$ & $2.45\times10^{-2}$& \boldmath{$5.20 \times 10^{-3}$} \\
       & 20 & \boldmath{$2.24\times10^{-2}$} &\boldmath{$2.14\times10^{-2}$} & $6.90 \times 10^{-3}$\\
    \bottomrule
    \end{tabular}
    \caption{RMSE of trained NDDEs of different sizes on testing data for the delay logistic equation~\eqref{eq:logistic}. 
    All NDDEs have two hidden layers with the same number of neurons in each layer as indicated in the second column.}
    \label{tab:logistic_NDDE}
\end{table}

\begin{table}[!t]
\small
    \centering    
    \begin{tabular}{llrrrrr}
    \toprule
   Loss & Size  & RMSE$_{x^\prime}$ & RMSE$_{x}$ & $|\tau-\hat{\tau}|$& Iter. & Time (s)\\ 
       \midrule
      \multirow{4}{*}{$\mathcal{L}_{x^\prime}$}  &  2 &  $1.63\times10^{-1}$ & $8.13\times10^{-2}$ & $1.14 \times 10^{-2}$ & \multirow{4}{*}{2000} & \multirow{4}{*}{$\approx$ 5}\\
       & 5 & $4.13\times10^{-2}$ & $1.14\times10^{-2}$ & $1.70 \times 10^{-3}$ && \\ 
       & 10 & $3.40\times10^{-2}$  & $7.90\times10^{-3}$  & $2.40 \times 10^{-3}$ & & \\ 
       & 20 & \boldmath{$3.07\times10^{-2}$} & \boldmath{$7.70\times10^{-3}$} & \boldmath{$1.60 \times 10^{-3}$} & & \\ 
     \midrule
     \multirow{4}{*}{$\mathcal{L}_{x}$}  &  2 &  $1.36\times10^{-1}$  & $5.24\times10^{-2}$ & $5.30 \times 10^{-3}$& \multirow{4}{*}{2000} & \multirow{4}{*}{$\approx$ 100}\\ 
       & 5 & $3.53\times10^{-2}$ & $1.17\times10^{-2}$ & $5.10 \times 10^{-3}$ & & \\
       & 10 & $3.37\times10^{-2}$ & $1.78\times10^{-2}$ & \boldmath{$2.30 \times 10^{-3}$} & & \\ 
       & 20 & \boldmath{$2.31\times10^{-2}$} & \boldmath{$3.70\times10^{-3}$} & $2.80 \times 10^{-3}$  &  & \\ 
    \bottomrule
    \end{tabular}
    \caption{RMSE of trained NDDEs of different sizes on testing data for the Mackey-Glass equation~\eqref{eq:mg}. 
    All NDDEs have two hidden layers with the same number of neurons in each layer as indicated in the second column.}
    \label{tab:mg_NDDE}
\end{table}

\begin{table}[!t]
\small
    \centering
    \begin{tabular}{cccccc}
    \toprule
    Loss & Size  & RMSE$_{x^\prime}$ & RMSE$_{x}$ & $|\tau_1-\hat{\tau}_1|$  &$|\tau_2-\hat{\tau}_2|$ \\
     &   & &  & ($\tau_1 = 0.0958$) & ($\tau_2 = 0.4792$)\\
       \midrule
       \multirow{4}{*}{$\mathcal{L}_{x^\prime}$} & 2 & $4.94\times10^{-1}$ & $4.81\times10^{-1}$  & $9.58\times10^{-2}$ & $1.17\times10^{-1}$ \\ 
     & 5 & $5.77\times10^{-2}$  & $1.93\times10^{-1}$  & $2.17\times10^{-2}$ & $7.50\times10^{-3}$\\
       & 10 & $2.93\times10^{-2}$ & $3.89\times10^{-2}$  & \boldmath{$1.58\times10^{-2}$} &\boldmath{$4.00\times10^{-3}$}\\
       & 20 & \boldmath{$1.77\times10^{-2}$} & \boldmath{$1.55\times10^{-2}$} & $4.18\times10^{-2}$ & $5.90\times10^{-3}$ \\ 
    \bottomrule
    \end{tabular}
    \caption{RMSE of trained NDDEs of different sizes on testing data for the climate model \eqref{eq:climate_model}. 
    All NDDEs have two hidden layers with the same number of neurons in each layer as indicated in the second column.}
    \label{tab:climate_NDDE}
\end{table}

For the Mackey-Glass example, we also compare the results for different loss functions in Table~\ref{tab:mg_NDDE}.
Once the delay is correctly learned, the performance of the trained networks are similar, regardless of the loss function. 
The computation time is also examined in this case.
Since the size of the networks is still relatively small, even for 20 neurons per layer, the main factor determining the computation time are the batchsize $N$ and the number of iterations.
For the simulation loss function, it also depends on the simulation horizon $H$. 
Observe that using derivative loss leads to a much faster training process compared to using simulation loss.

Finally, we remark that the generalizability of NDDEs was studied in~\citep{ji2024nd} through bifurcation diagrams and by computing the  maximum Lyapunov exponents. 
It was shown that using data with rich dynamics helps the generalizability, and that decoupling the delay and the nonlinearities allows the data-driven model to perform well for different delay values without re-training.

\CCLsection{Comparison of SINDy and NDDE approaches}
\label{sec:compare}

In this section, we evaluate the SINDy approach and the NDDE approach on the 3-state delayed R\"ossler system. 
The reconstructions of the time delays and the trajectories, as well as the computation time, are compared between different models and training methods.
In addition, the advantages and shortcomings of the different approaches are also discussed in the section.

\begin{example}
The double delayed \textit{R\"ossler system}~\citep{wu23,ghosh2008multiple,xu2010bifurcation} is given by
\begin{equation}\label{eq:rossler}
\left\{\setlength\arraycolsep{0.1em}\begin{array}{rcl}
x_1'(t) &=& -x_2(t) - x_3(t) + \alpha_1 x_1(t-\tau_1) + \alpha_2 x_1(t-\tau_2), 
\\[1mm]
x_2'(t) &=& x_1(t) + \beta_1 x_2(t), 
\\[1mm]
x_3'(t) &=& \beta_2 + x_3(t) x_1(t) - \gamma x_3(t).
\end{array}\right.
\end{equation}
This is a variation of the R\"ossler system, which simplifies the Lorenz model of turbulence. 
When ${\alpha_1 = 0.2}$, ${\alpha_2 = 1}$, ${\beta_1 = \beta_2 = 0.2}$, ${\gamma = 1.2}$ and ${\tau_1 = 1}$, ${\tau_2 = 2}$, the system exhibits chaotic behavior.
\end{example}

We collect data generated when starting from the constant initial function ${\varphi(s)\equiv [\,1.5,0.4,0.9\,]^{\top}}$.
The resulting dataset is a single trajectory with a section of length ${T_{\mathrm{tr}}=30}$ used for training and a section of length ${T_{\mathrm{ts}}=5}$ used for testing. 
The data are sampled with ${\Delta t = 0.05}$ and central difference is used to calculate the state derivatives numerically.

\CCLsubsubsection{SINDy for DDE approach}

We optimize both delays $\tau_1$ and $\tau_2$ and compare the results when using different degrees of polynomial libraries with ${d = 2, 3, 4}$.
All of these libraries include the true second-degree nonlinearity that appears in the RHS of \eqref{eq:rossler}, but for ${d = 3}$ or ${d = 4}$ the library also contains many other terms that are missing from \eqref{eq:rossler}.
When using {E-SINDy}, all delays are simultaneously optimized, and when using {P-SINDy}, only the maximum delay is optimized. 
For {E-SINDy}, candidate delays for $\tau_1$ and $\tau_2$ are uniformly sampled from the intervals ${[0.1, 1.5]}$ and ${[1, 3]}$. 
For {P-SINDy}, the candidate for the maximum delay $\bar{\tau}$ is uniformly sampled from the interval ${[0.1, 2.5]}$. 
PS optimization is utilized to search for the unknown delays. 
For {P-SINDy}, the degree of pseudospectral collocation for the discretization from DDE to ODE is set to ${M=5}$.

\CCLsubsubsection{NDDE approach}

We consider the NDDE
\begin{equation}
\label{eq:NDDE_rosslerr}
\hat{x}^\prime(t) = 
  W_3 \tanh{\Big(W_2 \tanh{\big(W_1 \hat{\mathbf{x}}_t +b_1\big)+b_2\Big)}},
\end{equation}
with different inputs and different number of neurons in the hidden layers.
While training with derivative loss, the input vector is from data ($\mathbf{x}_t$ instead of $\hat{\mathbf{x}}_t$ in \eqref{eq:NDDE_rosslerr}) and the output vector is the predicted state derivative $\hat{x}^\prime(t)$.
Depending on whether we know that the delays only enter the state $x_1$, or consider that they may enter $x_1$, $x_2$ and $x_3$, we can construct different network inputs: 
\begin{equation}\label{eq:simplifiedinput}
\mathbf{x}_t=
\begin{bmatrix}
      x_1(t) \\ x_2(t) \\ x_3(t) \\ {x_1(t-\hat{\tau}_1)} \\ {x_1(t-\hat{\tau}_2)}
  \end{bmatrix},
\qquad 
\mathbf{x}_t=
\begin{bmatrix}
      x_1(t) \\ x_2(t) \\ x_3(t) \\
  x_1(t-\hat{\tau}_1) \\ x_2(t-\hat{\tau}_1) \\ x_3(t-\hat{\tau}_1) \\
  x_1(t-\hat{\tau}_2) \\ x_2(t-\hat{\tau}_2) \\ x_3(t-\hat{\tau}_2)
  \end{bmatrix}.
\end{equation}
We refer the first input as simplified input and the second as full input and we will use both networks to learn from training data considering the derivative loss \eqref{eq:loss_deriv}.

\CCLsubsubsection{Comparison of the approaches}


\begin{figure}[!t]
\centering
\includegraphics[width=1\linewidth]{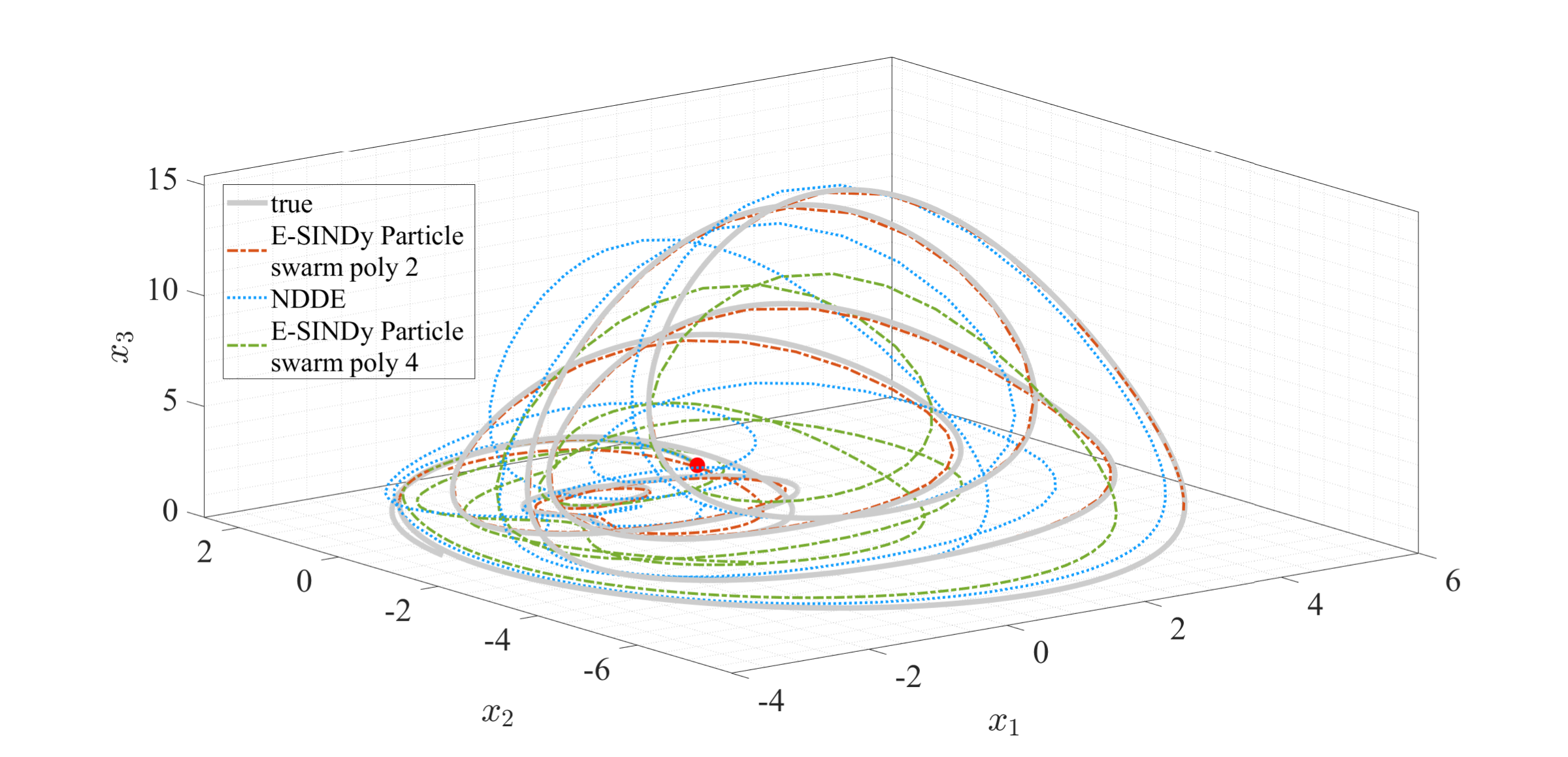}
\caption{Trajectory reconstruction for the R\"ossler model \eqref{eq:rossler} showing the true trajectory (gray), {E-SINDy} trajectory when using library with polynomial degree 2 (red),
{E-SINDy} trajectory when using library with polynomial degree 4 (green), and NDDE trajectory (blue). 
The true delay values are ${(\tau_1,\tau_2)=(1,2)}$. The optimized delay values obtained through particle swarm optimization for {E-SINDy} are ${(\widehat{\tau}_{1},\widehat{\tau}_{2})=(1.0001,2.0000)}$ when using polynomial library of order 2 and ${(\widehat{\tau}_{1},\widehat{\tau}_{2})=(1.0004,2.0000)}$ when using polynomial library of order 4. 
The NDDE approach yields ${(\widehat{\tau}_{1},\widehat{\tau}_{2})=(0.2309,1.9064)}$ when using simplified input and 5 neurons in each hidden layer. }\label{fig:Rossler-two-delay}
\end{figure}

\begin{figure}[!t]
\centering
\includegraphics[width=0.9\linewidth]{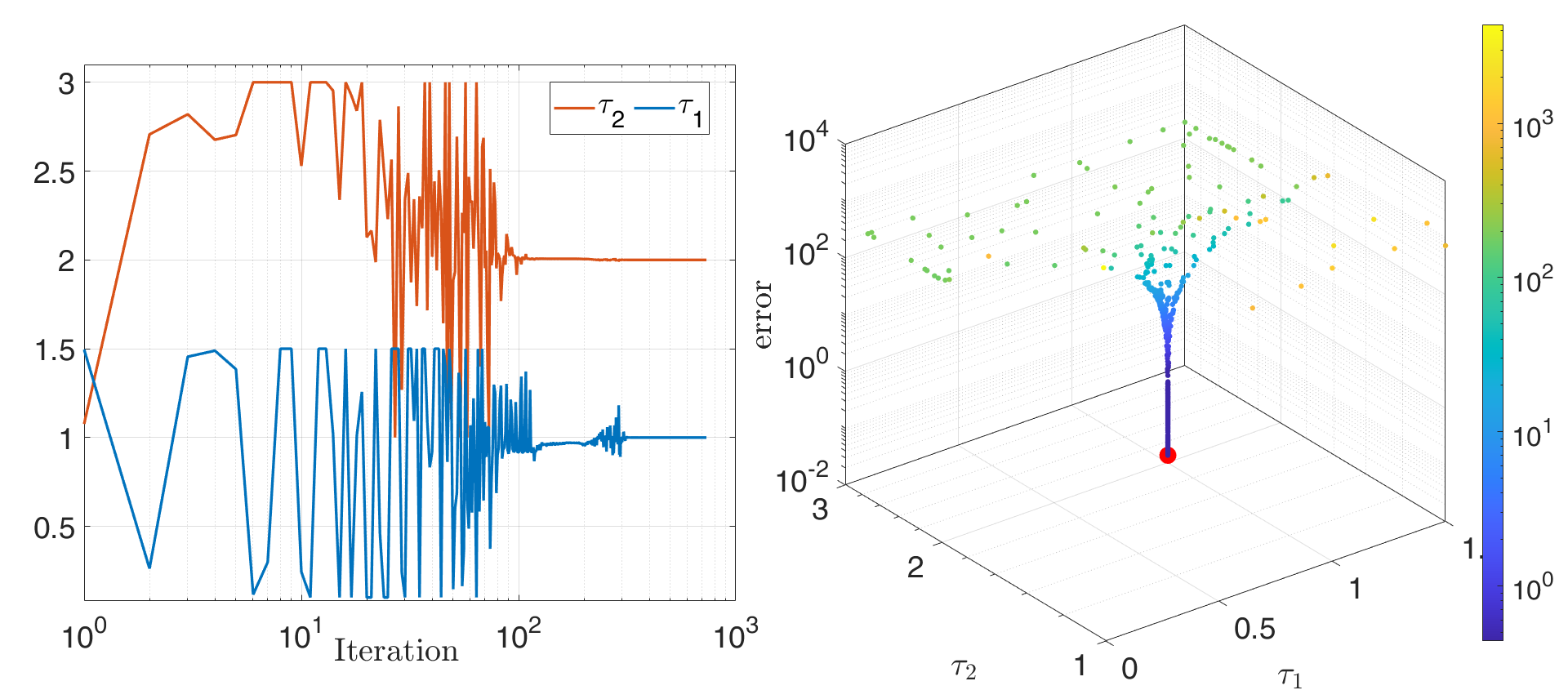}
\caption{Learned delays via particle swarm and particle filter errors along the iterations when applying {E-SINDy} with polynomial library of order 2 to the R\"ossler model \eqref{eq:rossler}.}\label{fig:rossler_delay_opt} 
\end{figure}

\begin{figure}[!t]
    \centering 
    \includegraphics[width=0.8\linewidth]{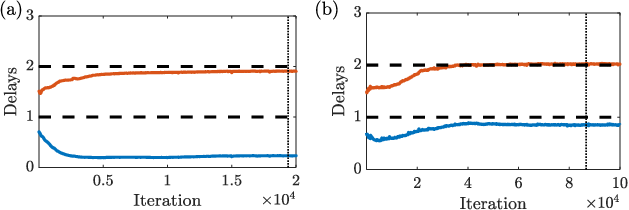}
    \caption{Learned delays along the iterations for the two-delay case when applying NDDE \eqref{eq:NDDE_rosslerr} to the R\"ossler model \eqref{eq:rossler}. (a) Simplified input and 5 neurons in each hidden layer, trained on one trajectory. (b) Simplified input and 20 neurons in each hidden layer, trained on ten trajectories.}\label{fig:rossler_2delay_ndde}
\end{figure}

The trajectory recovery is shown in Figure~\ref{fig:Rossler-two-delay} for the {E-SINDy} approach (with PS and polynomial libraries of orders 2 and 4) and for the NDDE approach (with simplified input, 5 neurons per layer and derivative loss). 
When the degree of the SINDy library is 2, it outperforms the NDDE in terms of prediction accuracy, while for degree 4 the SINDy predictions are less accurate than the NDDE predictions. 

The {E-SINDy} approach with low order library can recover the data very well since the components at the RHS of \eqref{eq:rossler} exist in the library and no higher-order nonlinearities are attempted to fit the data. 
However, when these higher-order terms are also included, the prediction accuracy decreases, as {E-SINDy} tries to capture the behavior by blending polynomials of different degrees.
In the meantime, the PS optimization still accurately captures the time delays in both cases.
For the limited amount of data utilized the accuracy of NDDEs predictions fall between the two SINDy cases while the learned delays remain somewhat far from the true values.
However, by including more training data, the performance of the NDDE can be significantly improved as will be demonstrated below. 

The training procedure for {E-SINDy} with particle swarm optimization is shown in Figure~\ref{fig:rossler_delay_opt}.
One may observe a clear minimum around the true delay values, which can be approached within a few hundred iterations
with the help of the PS optimizer.
For the NDDE, we depict the iterations in the trainable delays for different training datasets in Figure~\ref{fig:rossler_2delay_ndde}.
While for the limited data described above only one of the two delays is captured accurately, see panel (a), both delays can be recovered when the amount of training data is increased 10 times (by using 10 different initial conditions to generate 10 distinct trajectories), as shown in panel (b).

\begin{sidewaystable}[htbp]
\centering
\small
\begin{tabular}{lrrrrrrrrr}
\toprule
SINDy &  &  & {Training} &\multicolumn{2}{c}{Testing} & Calls &Time (s) \\
\cmidrule(r){5-6}
&$|\tau_{1}-\widehat{\tau}_1|$ & $|\tau_{2}-\widehat{\tau}_2|$ & RMSE$_{x^\prime}$ & RMSE$_{x^\prime}$ & RMSE$_{x}$   &&  \\
\midrule
\textbf{Polynomial order 2} & && & & && \\
E & $6.73\times10^{-3}$ & $3.16\times10^{-5}$ &{$4.31\times10^{-3}$}& {$2.82\times10^{-3}$}  & $8.17\times10^{-2}$& 704 & \textbf{86} \\ 

P5 & - & \boldmath{$3.63\times10^{-7}$} & \boldmath{$5.33\times10^{-5}$}  & \boldmath{$6.49\times10^{-5}$} & \boldmath{$5.84\times10^{-4}$} &\textbf{297} & 156 \\
\midrule
\textbf{Polynomial order 3} & && & && & \\
E & $5.89\times10^{-4}$ & $9.36\times10^{-5}$ & {$1.55\times10^{-2}$} & {$4.88\times10^{-2}$} & \boldmath{$8.65\times10^{-1}$}  &603 & \textbf{191} \\

P5 & - & \boldmath{$2.75\times10^{-7}$} & \boldmath{$4.84\times10^{-3}$} & \boldmath{$6.38\times10^{-3}$} & $3.20\times10^{-2}$ & \textbf{332} & 304 \\
\midrule
\textbf{Polynomial order 4} & && & && & \\
E & $7.32\times10^{-4}$ & $1.04\times10^{-5}$ & \boldmath{$3.02\times10^{-1}$} & \boldmath{$2.64$} & {$5.36$}  &658 & {2143} \\

P5 & - & \boldmath{$6.42\times10^{-7}$} & $7.66\times10^{-1}$ & $9.63$ & \boldmath{$3.28$} & \textbf{247} & \textbf{568} \\
\bottomrule
\end{tabular}%
\caption{Performance of {E-SINDy} and {P-SINDy} with ${M=5}$ and different polynomial orders for the R\"ossler system \eqref{eq:rossler}, spelling out the absolute error between the true and optimized delays, the RMSEs for training and testing phases, the number of SINDy calls and the CPU time while employing particle swarm.}\label{tab:comparison_rossler}
\end{sidewaystable} 

\begin{sidewaystable}[htbp]
    \centering
    \small
    \begin{tabular}{lrrrrrrr}
    \toprule
   NDDE & \multicolumn{2}{c}{}  & Training  & \multicolumn{2}{c}{Testing}  & & \\
 \cmidrule(r){5-6} 
    &$|\tau_{1}-\hat{\tau}_1|$ & $|\tau_{2}-\hat{\tau}_2|$   & $\mathcal{L}_{x^\prime}$ & RMSE$_{x^\prime}$ & RMSE$_{x}$ & Iter. & Time (s)\\
       \midrule
Full (5)  & $8.70\times10^{-1}$ &  $1.10$ & $4.53\times10^{-5}$ & \boldmath{$4.67\times10^{-2}$} & $1.86\times10^{-1}$ &\multirow{4}{*}{20000} & \multirow{4}{*}{$\approx$ 70} \\
   Full (20)    & $4.61\times10^{-1}$ & $6.32\times10^{-1}$ &  $5.53\times10^{-6}$ &  $7.09\times10^{-1}$ & $2.79$ & & \\
   Simplified (5)  & $7.69\times10^{-1}$ &  {$9.36\times10^{-2}$} & $1.91\times10^{-4}$ & $6.84\times10^{-2}$ & \boldmath{$8.96\times10^{-2}$} & & \\
  Simplified (20)  & $3.24\times10^{-1}$ & $3.29\times10^{-1}$ &   \boldmath{$4.15\times10^{-6}$} & $1.69\times10^{-1}$ & $1.87\times10^{-1}$ & & \\
  \cmidrule(r){7-8} 
  Simplified$^*$ (20)  & \boldmath{$1.39\times10^{-1}$} & \boldmath{$2.30\times10^{-2}$} &   $1.73\times10^{-5}$ & $2.21\times10^{-1}$ & $1.24\times10^{-1}$ & 100000 & $\approx$ 350\\
    \bottomrule
    \end{tabular}%
    \caption{Performance of NDDE for the R\"ossler system \eqref{eq:rossler}, spelling out the absolute error between the true and optimized delays and the errors in training and testing phases. The last row (Simplified$^*$) is trained on a larger dataset, where 10 trajectories are generated from 10 different constant histories.}
    \label{tab:rossler_NDDE}
\end{sidewaystable}

Table~\ref{tab:comparison_rossler} summarizes the SINDy results for E and P5 with PS, including the absolute errors on both delays $\tau_1$ and $\tau_2$, as well as the reconstruction error $\textrm{RMSE}_{x'}$ for training and the errors $\textrm{RMSE}_{x'}$ and $\textrm{RMSE}_{x}$ for testing. 
Different libraries of polynomial orders are compared in the table.
Both methods achieve excellent reconstruction accuracy for polynomial order 2. 
{E-SINDy} shows superior computational efficiency but {P-SINDy} achieves slightly better delay estimation accuracy in most cases.
When the order of polynomials increases in the library, the delay estimation remains accurate, while the reconstruction accuracy drops. 
This illustrates that having the right set of nonlinear functions in the library is important for SINDy to accurately identify the underlying dynamics of time delay systems.

The NDDE results are collected in Table~\ref{tab:rossler_NDDE}. 
For each network input choice, we test two different network sizes, i.e., the number of neurons in each hidden layer is either 5 or 20, as indicated in the first column. 
The networks with simple input structure lead to smaller errors in the learned delays, that is, having restriction on the network input helps regulating the network.
In addition, increasing the number of neurons in each hidden layer from 5 to 20 yields better predictions on the testing data.
Training the NDDE on a dataset consisting of 10 different trajectories results in more accurate learned delays and simulation results compared to the NDDE trained on only one trajectory.
That is, while learning from a single trajectory is challenging for neural networks (without any knowledge about the nonlinearities), utilizing data with multiple trajectories stemming from different initial histories improves the learning performance. 
In summary, having richer training data may compensate the missing apriori knowledge about the form of the RHS.

\CCLsection{Conclusions}

In this chapter, we presented state-of-art data-driven methods for delay differential equations (DDE), including the DDE extension of the sparse identification of nonlinear dynamics (SINDy) and the deployment of neural delay differential equations (NDDE) with trainable delays.
The core idea in both presented methods is to approximate the right-hand side (RHS) of DDEs, with meaningful delay(s). 

We illustrated the SINDy method designed for ordinary differential equations (ODE) as extended to DDEs and focused on two approaches, the expert SINDy ({E-SINDy}) and the pragmatic SINDy ({P-SINDy}). 
The difference exists in whether we approximate the RHS of the DDE directly or approximate its ODE representation. 
In the {E-SINDy}, the library is constructed for the DDE directly and all the delays need to be optimized. 
In the {P-SINDy}, the library is constructed based on the ODE discretization, therefore only the maximum delay needs to be identified. 
{P-SINDy} can be very useful when the number of delays increases, however, it may not give back an interpretable DDE from the ODE representation when there are multiple delays. 
In both SINDy approaches, the delays and function parameter in the library can be optimized externally, multiple optimizers (brute force, Bayesian optimization, particle swarm optimization) were tested.

For the neural network approach, we provided an overview on using neural networks in continuous time for learning the DDEs. 
The neural networks considered in this chapter are equipped with trainable delays and they result in interpretable form of NDDEs. 
Different from the SINDy approaches which use external optimizers for learning the delays, the NDDE approach can learn the delays and nonlinear dynamics simultaneously via gradient-based algorithms. 
Multiple loss functions were introduced. 
When using derivative loss, we focused on the training algorithms without simulating the networks, which resulted in better computational efficiency.
When using the simulation loss, undesired local minima (numerical differentiation schemes) can be avoided and better prediction accuracy can be achieved, while paying the price of slower training.

Multiple test examples were provided in the chapter and the corresponding MATLAB implementation codes were also released to public. 
We revealed both the advantages and the fundamental limitations of the data-driven methods. 
From the results, we can see the importance of having prior knowledge (i.e., what type of nonlinearities and up to which order need to be included in the SINDy library) and having good quality of training data (that cover large enough part of the state space for NDDEs).
By expressing the delays explicitly and constructing the RHS of the underlying DDEs, the data-driven models given by the E-SINDy approach and by the NDDE approach are interpretable.
On the other hand, the connection between DDEs and ODEs serves as a foundation for the P-SINDy approach and for the neural ordinary differential equation (NODE) approach pursued by~\cite{ji2022NODE} (not included in this chapter for brevity).
The resulting data-driven models are in the form which enables their analysis using well-developed tools of nonlinear dynamics and control.

In general, the goal of this chapter is not to promote any specific data-driven method over the others, but to provide different views to modeling time delay systems from data and to display the potentials of different methods. 
The choice of data-driven methods is problem-dependent and it should be customized based on the objective of the tasks.
There are still many open questions in this broad topic. 
For example, how to determine the quantity and the quality of the training data which is sufficient to learn the dynamics of time delay systems, and how this is related to the amount of knowledge embedded in the data-driven model. 
Studies on stochastic and distributed delays can be further explored and methods can also be extended to neutral delay differential equations. 
It will be interesting to see applications on real-world data, even though the ground truth is normally unknown. 
The data-driven methods can also be improved, for instance, restricting the NDDE for specific training tasks (to employ domain knowledge) and equipping SINDy approaches with gradient-based delay learning.

\bibliographystyle{plainnat}
\bibliography{Contribs}

@article{bpk16,
  title={Discovering governing equations from data by sparse identification of nonlinear dynamical systems},
  author={Brunton, S. L. and Proctor, J. L. and Kutz, J. N.},
  journal={Proceedings of the National Academy of Sciences},
  volume={113},
  number={15},
  pages={3932--3937},
  year={2016},
  publisher={National Academy Sciences}
}

@book{bk19,
  title={Data-driven Science and Engineering},
  author={Brunton, S. L. and Kutz, J. N.},
  year={2022},
  publisher={Cambridge University Press}
}

@article{rud17,
  title={Data-driven discovery of partial differential equations},
  author={Rudy, S. H. and Brunton, S. L. and Proctor, J. L. and Kutz, J. N.},
  journal={Science Advances},
  volume={3},
  number={4},
  pages={e1602614},
  year={2017},
  publisher={American Association for the Advancement of Science}
}

@article{bnc18,
  title={Sparse learning of stochastic dynamical equations},
  author={Boninsegna, L. and N{\"u}ske, F. and Clementi, C.},
  journal={The Journal of chemical physics},
  volume={148},
  number={24},
  year={2018},
  publisher={AIP Publishing}
}

@article{champ19,
  title={Data-driven discovery of coordinates and governing equations},
  author={Champion, K. and Lusch, B. and Kutz, J. N. and Brunton, S. L.},
  journal={Proceedings of the National Academy of Sciences},
  volume={116},
  number={45},
  pages={22445--22451},
  year={2019},
  publisher={National Acadademy Sciences}
}

@article{bary04,
  title={Barycentric lagrange interpolation},
  author={Berrut, J. P. and Trefethen, L. N.},
  journal={SIAM review},
  volume={46},
  number={3},
  pages={501--517},
  year={2004},
  publisher={SIAM}
}

@book{trefethen2000,
  title={Spectral methods in MATLAB},
  author={Trefethen, L. N.},
  year={2000},
  publisher={SIAM}
}

@article{chtd11,
  title={Numerical differentiation of noisy, nonsmooth data},
  author={Chartrand, R.},
  journal={International Scholarly Research Notices},
  volume={2011},
  number={1},
  pages={164564},
  year={2011},
  publisher={Wiley Online Library}
}

@article{lasso96,
  title={Regression shrinkage and selection via the lasso},
  author={Tibshirani, R.},
  journal={Journal of the Royal Statistical Society Series B},
  volume={58},
  number={1},
  pages={267--288},
  year={1996},
  publisher={Oxford University Press}
}

@article{sandoz23,
  title={SINDy for delay-differential equations: application to model bacterial zinc response},
  author={Sandoz, A. and Ducret, V. and Gottwald, G. A and Vilmart, G. and Perron, K.},
  journal={Proceedings of the Royal Society A},
  volume={479},
  number={2269},
  pages={20220556},
  year={2023},
  publisher={The Royal Society}
}

@book{krasovskii63,
  title={Problems of the Theory of Stability of Motion},
  author={Krasovskii, N.},
  year={1963},
  publisher={English translation, Stanford University Press}
}

@article{shayer2000,
  title={Stability, bifurcation, and multistability in a system of two coupled neurons with multiple time delays},
  author={Shayer, L. P. and Campbell, S. A.},
  journal={SIAM Journal on Applied Mathematics},
  volume={61},
  number={2},
  pages={673--700},
  year={2000},
  publisher={SIAM}
}

@article{ghosh2008multiple,
  title={Multiple delay {R}{\"o}ssler system—Bifurcation and chaos control},
  author={Ghosh, D. and Chowdhury, A. R. and Saha, P.},
  journal={Chaos, Solitons \& Fractals},
  volume={35},
  number={3},
  pages={472--485},
  year={2008},
  publisher={Elsevier}
}

@article{xu2010bifurcation,
  title={Bifurcation analysis of {R}{\"o}ssler system with multiple delayed feedback},
  author={Xu, M. and Wei, Y. and Wei, J.},
  journal={Electronic Journal of Qualitative Theory of Differential Equations},
  volume={2010},
  number={63},
  pages={1--22},
  year={2010},
  publisher={University of Szeged, Hungary}
}

@book{diekmann95,
  title={Delay equations: functional-, complex-, and nonlinear analysis},
  author={Diekmann, O. and Van Gils, S. A. and Lunel, S. M. V. and Walther, H. O.},
  volume={110},
  year={1995},
  publisher={Applied Mathematical Sciences. Springer Verlag}
}

@article{pec24,
  title={Data-driven discovery of delay differential equations with discrete delays},
  author={Pecile, A. and Demo, N. and Tezzele, M. and Rozza, G. and Breda, D.},
  journal={Journal of Computational and Applied Mathematics},
  volume={461},
  pages={116439},
  year={2025},
  publisher={Elsevier}
}

@article{bdgsv16,
  title={Pseudospectral discretization of nonlinear delay equations: new prospects for numerical bifurcation analysis},
  author={Breda, D. and Diekmann, O. and Gyllenberg, M. and Scarabel, F. and Vermiglio, R.},
  journal={SIAM Journal on applied dynamical systems},
  volume={15},
  number={1},
  pages={1--23},
  year={2016},
  publisher={SIAM}
}

@inproceedings{kopeczi23,
  title={Data-driven delay identification with SINDy},
  author={K{\"o}peczi-B{\'o}cz, Tam{\'a}s, A. and Sykora, H. and Tak{\'a}cs, D.},
  booktitle={International Conference on Nonlinear Dynamics and Applications},
  pages={481--491},
  year={2023},
  organization={Springer}
}

@article{borgioli2020pseudospectral,
  title={Pseudospectral method for assessing stability robustness for linear time-periodic delayed dynamical systems},
  author={Borgioli, F. and Hajdu, D. and Insperger, T. and Stepan, G. and Michiels, W.},
  journal={International Journal for Numerical Methods in Engineering},
  volume={121},
  number={16},
  pages={3505--3528},
  year={2020},
  publisher={Wiley Online Library}
}

@article{lehotzky2016pseudospectral,
  title={A pseudospectral tau approximation for time delay systems and its comparison with other weighted-residual-type methods},
  author={Lehotzky, D. and Insperger, T.},
  journal={International Journal for Numerical Methods in Engineering},
  volume={108},
  number={6},
  pages={588--613},
  year={2016},
  publisher={Wiley Online Library}
}

@article{Hutchinson1948,
	author = {Hutchinson, G. E.},
	journal = {Annals of the New York Academy of Sciences},
	number = {4},
	pages = {221-246},
	title = {Circular Causal Systems in Ecology},
	volume = {50},
	year = {1948}
}

@book{Gopalsamy1992,
	author = {Gopalsamy, K.},
	publisher = {Springer},
	title = {Stability and Oscillations in Delay Differential Equations of Population Dynamics},
	year = 1992
}

@article{bdls16,
  title={Numerical bifurcation analysis of a class of nonlinear renewal equations},
  author={Breda, D. and Diekmann, O. and Liessi, D. and Scarabel, F.},
  journal={Electronic Journal of Qualitative Theory of Differential Equations},
  volume={2016},
  number={65},
  pages={1--24},
  year={2016},
  publisher={University of Szeged, Hungary}
}

@article{spjc07,
  title={A simple chaotic delay differential equation},
  author={Sprott, J. C.},
  journal={Physics Letters A},
  volume={366},
  number={4-5},
  pages={397--402},
  year={2007},
  publisher={Elsevier}
}

@book{hale77,
	author = {{Hale}, J. K. and {Verduyn Lunel}, S. M.},
	publisher = {Springer},
	title = {Introduction to Functional Differential Equations},
	year = {1993}
}

@article{bbt24,
  title={Sparse identification of time delay systems via pseudospectral collocation},
  author={Bozzo, E. and Breda, D. and Tanveer, M. },
  journal={IFAC-PapersOnLine},
  volume={58},
  number={27},
  pages={108--113},
  year={2024},
  publisher={Elsevier}
}

@article{mg77,
  title={Oscillation and chaos in physiological control systems},
  author={Mackey, M. C. and Glass, L.},
  journal={Science},
  volume={197},
  number={4300},
  pages={287--289},
  year={1977},
  publisher={American Association for the Advancement of Science}
}

@article{wu23,
  title={Reconstruction of delay differential equations via learning parameterized dictionary},
  author={Wu, Y.},
  journal={Physica D},
  volume={446},
  pages={133647},
  year={2023},
  publisher={Elsevier}
}

@book{garnett2023,
  author    = {Garnett, R.},
  title     = {{Bayesian Optimization}},
  year      = {2023},
  publisher = {Cambridge University Press}
}

@book{smith2011introduction,
  title={An Introduction to Delay Differential Equations with Applications to the Life Sciences},
  author={Smith, H.},
  year={2011},
  publisher={Springer}
}

@article{shahriari2015,
  title={Taking the human out of the loop: A review of Bayesian optimization},
  author={Shahriari, B. and Swersky, K. and Wang, Z. and Adams, R. P. and De Freitas, N.},
  journal={Proceedings of the IEEE},
  volume={104},
  number={1},
  pages={148--175},
  year={2015},
  publisher={IEEE}
}

@book{williams2006,
  title={Gaussian Processes for Machine Learning},
  author={Williams, C. K. I. and Rasmussen, C. E.},
  year={2005},
  publisher={MIT Press}
}

@article{bonyadi2017,
  title={Particle swarm optimization for single objective continuous space problems: a review},
  author={Bonyadi, M. R. and Michalewicz, Z.},
  journal={Evolutionary Computation},
  volume={25},
  number={1},
  pages={1--54},
  year={2017},
  publisher={MIT Press}
}

@inproceedings{kennedy1995,
  title={Particle swarm optimization},
  author={Kennedy, J. and Eberhart, R.},
  booktitle={International Conference on Neural Networks},
  volume={4},
  pages={1942--1948},
  year={1995},
  organization={IEEE}
}

@book{erneux2009applied,
  title={Applied Delay Differential Equations},
  author={Erneux, T.},
  publisher={Springer},
  year={2009}
}

@inproceedings{shi1998,
  title={A modified particle swarm optimizer},
  author={Shi, Y. and Eberhart, R.},
  booktitle={IEEE International Conference on Evolutionary Computation},
  pages={69--73},
  year={1998},
  organization={IEEE}
}

@inproceedings{brute2011,
  title={A brute-force algorithm for reconstructing a scene from two projections},
  author={Enqvist, O. and Jiang, F. and Kahl, F.},
  booktitle={CVPR 2011},
  pages={2961--2968},
  year={2011},
  organization={IEEE}
}

@mastersthesis{brute2019,
  title={The Brute Force Algorithm},
  author={Weyer, A.},
  year={2019},
  school={Bowling Green State University}
}

@book{engel2000one,
  title={One-parameter semigroups for linear evolution equations},
  author={Engel, K.-J. and Nagel, R. and Brendle, S.},
  year={2000},
  publisher={Springer}
}

@book{stepan1989retarded,
  title={Retarded Dynamical Systems: Stability and Characteristic Functions},
  author={St{\'e}p{\'a}n, G.},
  publisher={Longman},
  year={1989}
}

@book{insperger2011semi,
  title={Semi-discretization for Time-delay Systems: Stability and Engineering Applications},
  author={Insperger, T. and St{\'e}p{\'a}n, G.},
  year={2011},
  publisher={Springer}
}

@book{breda2022controlling,
  title={Controlling Delayed Dynamics: Advances in Theory, Methods and Applications},
  author={Breda, D.},
  volume={604},
  year={2022},
  publisher={Springer}
}

@article{bellen2009recent,
  title={Recent trends in the numerical solution of retarded functional differential equations},
  author={Bellen, A. and Maset, S. and Zennaro, M. and Guglielmi, N.},
  journal={Acta Numerica},
  volume={18},
  pages={1--110},
  year={2009},
  publisher={Cambridge University Press}
}

@book{kuznetsov1998elements,
  title={Elements of Applied Bifurcation Theory},
  author={Kuznetsov, Y. A.},
  year={1998},
  publisher={Springer}
}

@Article{pei2013mapping,
  author    = {Pei, J.-S. and Mai, E. C. and Wright, J. P. and Masri, S. F.},
  journal   = {Nonlinear Dynamics},
  title     = {Mapping some basic functions and operations to multilayer feedforward neural networks for modeling nonlinear dynamical systems and beyond},
  year      = {2013},
  number    = {1-2},
  pages     = {371--399},
  volume    = {71},
  publisher = {Springer},
}

@article{canziani2016analysis,
  title={An analysis of deep neural network models for practical applications},
  author={Canziani, A. and Paszke, A. and Culurciello, E.},
  journal={arXiv preprint arXiv:1605.07678},
  year={2016}
}

@InProceedings{ji2020, 
title = {Feed-forward Neural Networks with Trainable Delay}, 
author = {Ji, X. A. and Moln\'ar, T. G. and Avedisov, S. S. and Orosz, G.}, 
booktitle = 	{2nd Conference on Learning for Dynamics and Control},
pages = {127--136}, 
year = {2020}, 
volume = {120}, 
publisher =   {PMLR}
}

@article{levine2022framework,
  title={A framework for machine learning of model error in dynamical systems},
  author={Levine, M. and Stuart, A.},
  journal={Communications of the American Mathematical Society},
  volume={2},
  number={07},
  pages={283--344},
  year={2022}
}

@article{gupta2023neural,
  title={Generalized neural closure models with interpretability},
  author={Gupta, A. and Lermusiaux, P. F. J.},
  journal={Scientific Reports},
  volume={13},
  pages={10634},
  year={2023},
}

@article{keane2019effect,
  title={The effect of state dependence in a delay differential equation model for the El Ni{\~n}o Southern Oscillation},
  author={Keane, A. and Krauskopf, B. and Dijkstra, H. A.},
  journal={Philosophical Transactions of the Royal Society A},
  volume={377},
  number={2153},
  pages={20180121},
  year={2019},
  publisher={The Royal Society Publishing}
}

@article{keane2017climate,
  title={Climate models with delay differential equations},
  author={Keane, A. and Krauskopf, B. and Postlethwaite, C. M.},
  journal={Chaos},
  volume={27},
  number={11},
  pages={114309},
  year={2017},
  publisher={AIP Publishing LLC}
}

@article{saunders2001boolean,
  title={A Boolean delay equation model of ENSO variability},
  author={Saunders, A. and Ghil, M.},
  journal={Physica D},
  volume={160},
  number={1-2},
  pages={54--78},
  year={2001},
  publisher={Elsevier}
}

@article{ji2022NODE,
  title={Learning Time Delay Systems with Neural Ordinary Differential Equations},
  author={Ji, X. A. and Orosz, G.},
  journal={IFAC-PapersOnLine},
  volume={55},
  number={36},
  pages={79--84},
  year={2022},
  publisher={Elsevier}
}

@inproceedings{ji2023l4dc,
  title={Learning the dynamics of autonomous nonlinear delay systems},
  author={Ji, X. A. and Orosz, G.},
  booktitle={3rd Conference on Learning for Dynamics and Control Conference},
  pages={116--127},
  year={2023},
  volume = {211}, 
  publisher={PMLR}
}

@inproceedings{ji2021l4dc,
  title={Learning the dynamics of time delay systems with trainable delays},
  author={Ji, X. A. and Moln{\'a}r, T. G. and Avedisov, S. S. and Orosz, G.},
  booktitle={Learning for Dynamics and Control},
  pages={930--942},
  volume = {144}, 
  year={2021},
  organization={PMLR}
}

@article{ji2024nd,
  title={Learn from one and predict all: single trajectory learning for time delay systems},
  author={Ji, X. A. and Orosz, G.},
  journal={Nonlinear Dynamics},
  volume={112},
  number={5},
  pages={3505--3518},
  year={2024},
  publisher={Springer}
}

@inproceedings{glorot2010understanding,
  title={Understanding the difficulty of training deep feedforward neural networks},
  author={Glorot, X. and Bengio, Y.},
  booktitle={13th International Conference on Artificial Intelligence and Statistics},
  pages={249--256},
  year={2010},
  volume={9},
  organization={PMLR}
}

@inproceedings{kingma2014adam,
  author    = {D. P. Kingma and
               J. Ba},
  editor    = {Y. Bengio and
               Y. LeCun},
  title     = {Adam: {A} Method for Stochastic Optimization},
  booktitle = {3rd International Conference on Learning Representations},
  year      = {2015},
}
\end{document}